\documentclass{amsart}

\usepackage{amsmath}
\usepackage{amssymb}
\usepackage{mathrsfs}

\newtheorem{theorem}{Theorem}[section]
\newtheorem{claim}[theorem]{Claim}

\newtheorem{conclusion}[theorem]{Conclusion}
\newtheorem{observation}[theorem]{Observation}

\theoremstyle{definition}
\newtheorem{definition}[theorem]{Definition}

\newtheorem{discussion}[theorem]{Discussion}

\newtheorem{hypothesis}[theorem]{Hypothesis}

\theoremstyle{remark}
\newtheorem{remark}[theorem]{Remark}
\newtheorem{question}[theorem]{Question}
\newtheorem{notation}[theorem]{Notation}

\newcommand{\otp}{{\rm otp}}
\newcommand{\lub}{{\rm lub}}

\newcommand{\st}{{\rm st}}

\newcommand{\gn}{{\rm gn}}

\newcommand{\Reg}{{\rm Reg}}

\newcommand{\rang}{{\rm rang}}

\newcommand{\Ax}{{\rm Ax}}

\newcommand{\COM}{{\rm COM}}
\newcommand{\INC}{{\rm INC}}

\newcommand{\Min}{{\rm Min}}
\newcommand{\Dom}{{\rm Dom}}
\newcommand{\Rang}{{\rm Rang}}
\newcommand{\stt}{{\rm st}}
\newcommand{\snn}{{\sn}}

\newcommand{\eqdf}{{=^{\rm df}}}
\newcommand{\uhr}{{\upharpoonright}}
\newcommand{\hz}{{}}
\newcommand{\nz}{{}}

\newcommand{\hh}{{\mathbf{h}}}
\newcommand{\hp}{{\mathbf{p}}}
\newcommand{\hn}{{\mathbf{\varepsilon }}}
\newcommand{\hm}{{\mathbf{\iota}}}
\newcommand{\hc}{{\mathbf{c}}}
\newcommand{\hi}{{\alpha}}

\newcommand{\rest}{{\restriction}}
\newcommand{\dom}{{\rm dom}}

\newcommand{\wilog}{{\rm without loss of generality}}

\newcommand{\then}{{\underline{then}}}
\newcommand{\when}{{\underline{when}}}

\newcommand{\Then}{{\underline{Then}}}

\newcommand{\Iff}{{\underline{iff}}}
\newcommand{\mn}{{\medskip\noindent}}
\newcommand{\sn}{{\smallskip\noindent}}

\newcommand{\bbR}{{\mathbb R}}

\newcommand{\cI}{{\mathscr I}}

\newcommand{\bbP}{{\mathbb P}}

\newcommand{\bbQ}{{\mathbb Q}}

\newcommand{\cS}{{\mathscr S}}

\newcommand{\cW}{{\mathscr W}}

\newcommand{\cU}{{\mathscr U}}

\newcommand{\cf}{{\rm cf}}

\newcount\skewfactor
\def\mathunderaccent#1#2 {\let\theaccent#1\skewfactor#2
\mathpalette\putaccentunder}
\def\putaccentunder#1#2{\oalign{$#1#2$\crcr\hidewidth
\vbox to.2ex{\hbox{$#1\skew\skewfactor\theaccent{}$}\vss}\hidewidth}}
\def\name{\mathunderaccent\tilde-3 }

\newenvironment{PROOF}[2][\proofname.]
   {\begin{proof}[#1]\renewcommand{\qedsymbol}{\textsquare$_{\rm #2}$}}
   {\end{proof}}

\usepackage[hidelinks]{hyperref}
\begin{document}
\makeatletter\def\shfiuwefootnote{\gdef\@thefnmark{}\@footnotetext}\makeatother\shfiuwefootnote{Version 2020-07-29. See \url{https://shelah.logic.at/papers/1036/} for possible updates.}

\title {Forcing axioms for $\lambda$-complete $\mu^+$-c.c. \\
Sh1036}
\author {Saharon Shelah}
\address{Einstein Institute of Mathematics\\
Edmond J. Safra Campus, Givat Ram\\
The Hebrew University of Jerusalem\\
Jerusalem, 91904, Israel\\
 and \\
 Department of Mathematics\\
 Hill Center - Busch Campus \\
 Rutgers, The State University of New Jersey \\
 110 Frelinghuysen Road \\
 Piscataway, NJ 08854-8019 USA}
\email{shelah@math.huji.ac.il}
\urladdr{http://shelah.logic.at}
\thanks{Research supported by the US-Israel Binational Science Foundation
(Grant No. 2010405). References like \cite[Th0.2=Ly5]{Sh:950} means
the label of Th.0.2 is y5.  The reader should note that the version in
my website is usually more updated than the one in the mathematical archive.
The author thanks Alice Leonhardt for the beautiful typing.
  First typed July 31, 2012}




\date{2020-07-29}

\begin{abstract}
We consider forcing axioms for suitable families of
$\mu$-complete $\mu^+$-c.c. forcing notions.
We show that some form of the condition ``$p_1,p_2$ have a
$\le_{\bbQ}-\lub$ in $\bbQ$" is necessary.  We also show some
versions are really stronger than others.
\end{abstract}

\maketitle
\numberwithin{equation}{section}
\setcounter{section}{-1}
\newpage

\section{Introduction}
\bigskip

\subsection {Is Well Met Necessary in Some Forcing Axiom}\
\bigskip

We investigate the relationships between some forcing axioms related
to pressing down functions for $\mu^+$-c.c., mainly
from \cite{Sh:546}.  This in particular is to answer
Kolesnikov's question of having $\bbP$ satisfying one condition but
with  
no $\bbP'$
equivalent to $\bbP$ satisfying another.   A side issue is clarifying
a point in \cite{Sh:927}
(a rephrasing is $(2)^\varepsilon_{c,D}$ from \ref{x3}).  We intend to
continue this considering related axioms in \cite{Sh:F1856}.

We justify the ``well met, having lub" in some forcing axioms, e.g. condition
(c) in $*^1_{\mu,\bbQ}$.

In \cite{Sh:80} such forcing axiom was proved consistent, for forcing
notion satisfying (for $ \mu ^ {< \mu }= \mu $; 
we may write $ \text{`}\text{`})\mathbb{Q} $ satisfies 
$ * ^1_\mu $ instead $ * ^1_{\mu, \mathbb{Q} }$ , similarly below):
\mn
\begin{enumerate}
\item[$*^1_{\mu,\bbQ}$]  $\bbQ$
is 
a forcing notion such that:
\sn
\begin{enumerate}
\item  $(< \mu)$-complete, i.e. any
increasing sequence of length $<\mu$ has an upper bound
\sn
\item  $\mu^+$-regressive-c.c.: if $p_\alpha \in \bbQ$ for $\alpha
  <\mu^+$ then for some club $E$ of $\mu^+$ and pressing down function
  $f$ on $E$ we have $[\delta_1 \in E \wedge \delta_2 \in E \wedge 
  \big( f(\delta_1) =
  f(\delta_2) \big) 
  \wedge \big(  \cf( \delta _1) = \mu = \cf( \delta _2 ) \big)
  \Rightarrow p_{\delta_1},p_{\delta_2}$ compatible]
\sn
\item  if $p_1,p_2 \in \bbQ$ are compatible then $p_1,p_2$ have a lub.
\end{enumerate}
\end{enumerate}
\mn
An easily stated version which is still enough is:
\mn
\begin{enumerate}
\item[$*^2_{\mu,\bbQ}$]  $\bbQ$ is a forcing notion satisfying clause (a) and
\sn
\begin{enumerate}
\item[$(b)'$]  if $p_\alpha \in \bbQ$ for $\alpha < \mu^+$ \then \, for
some $(E,\bar q,f)$ we have
\sn
\item[${{}}$]  $\bullet \quad E$ a club of $\mu^+$
\sn
\item[${{}}$]  $\bullet \quad \bar q = \langle q_\alpha:\alpha <
\mu^+\rangle$
\sn
\item[${{}}$]  $\bullet \quad p_\alpha \le_{\bbQ} q_\alpha$
\sn
\item[${{}}$]  $\bullet \quad f$ is a pressing down function on $E$
\sn
\item[${{}}$]  $\bullet \quad$ if $\delta_1 \in E \wedge \delta_2 \in
E \wedge \cf(\delta_1) = \mu = \cf(\delta_2) \wedge f(\delta_1) =
f(\delta_2)$ \then \,

\hskip30pt  $q_{\delta_1},q_{\delta_2}$ has a lub.
\end{enumerate}
\end{enumerate}
\mn
An obvious fact used is
\mn
\begin{enumerate}
\item[$\boxplus$]  Assume $\bbQ$ is a
forcing notion, 
$\varepsilon < \mu$ a limit
  ordinal, $\bar p_\ell = \langle p_{\ell,\alpha}:\alpha < \varepsilon
\rangle$ is $\le_{\bbQ}$-increasing for $\ell=1,2$ and 
for every $ \alpha < \varepsilon $
the condition 
$p_\alpha
\in \bbQ$ is a $\le_{\bbQ}$-lub of $p_{1,\alpha},p_{2,\alpha}$
  (i.e. $\bigwedge\limits^2_{\ell=1} p_{\ell,\alpha} \le_{\bbQ}
p_\alpha$ and $(\forall q)(p_{1,\alpha} \le_{\bbQ} q
\wedge p_{2,\alpha} \le_{\bbQ} q \Rightarrow p_\alpha \le_{\bbQ} q))$.  
\Then \, $\langle p_\alpha:\alpha <
\varepsilon \rangle$ is $\le_{\bbQ}$-increasing, hence
if $\{p_\alpha:\alpha < \varepsilon\}$ has an upper bound then so does
  $\{p_{1,\alpha},p_{2,\alpha}:\alpha < \varepsilon\}$.
\end{enumerate}
\mn

Now \cite{Sh:963} 
mainly deal with consistency results  
for singular $ \mu $, but on the way   
has (with a complete proof of the iteration theorem)
suggest a condition weaker than the one in \cite{Sh:80} and even the one
in \cite{Sh:93}  and is stronger than the one in  
\cite[1.7(1)]{Sh:546}, using a trivial strategy and $ \varepsilon = \omega $.
Using \ref{x2}, 
the condition  from \cite{Sh:546} is $ (2)^\varepsilon _ {c, D}$,  
where $ \varepsilon $ is a limit  ordinal $ < \mu $, 
and the condition from \cite{Sh:963} is

\begin{enumerate} 
\item[$*^3_{\mu,\bbQ}$]  $\bbQ$ a forcing notion such that
\sn
\begin{enumerate}
\item[(a)] as above
\item[(b)] as above 
\item[(c)]  if, for every $ n < \omega $ we have
$ p_n \le p_{n + 1 }, q_n \le q _{n + 1 }$ and $ p_n , q_n $ are
compatible \then \, the set $ \{ p_n, q_n : n < \omega \} $ 
has a common upper bound (here this is clause $ (3)_{b,\omega }$ 
of Def \ref{x2}).
\end{enumerate} 
\end{enumerate}

Our main conclusions are \ref{a20}, \ref{a24}, \ref{c3}, \ref{d25}.

The immediate reason for this paper is that the statement in
Baldwin-Kolesnikov-Shelah \cite[3.6]{Sh:927} is
misquoting \cite[4.12]{Sh:93}.  
We shall show below that the statement
is inconsistent because as stated it totally waives the condition ``every two
compatible members of $\bbP$ have a lub".  Also, it is stated that in
\cite[4.12]{Sh:93} this was claimed, but quoting only \cite{Sh:80}.  
In
Shelah-Spinas \cite{Sh:1110} we consider another strengthening of
the axioms.

More fully, \cite[4.12]{Sh:93} omits the condition above, but demands the
existence of lub's of some pairs of conditions so that it holds
in the cases it is actually used.
So, in that case  
the proof of \cite{Sh:80} works, and see more in
\cite[Def.1.1]{Sh:546} which gives an even weaker condition called
$*^\varepsilon_\mu$.

Concerning $*^1_{\mu,\bbQ}$, the preservation of a related condition
was proved independently by Baumgartner, who instead of (b) use 
somewhat stronger  condition $(b)^+$ which says that
$\bbQ$ is the union of $\mu$ sets of pairwise
compatible elements with lub, this  
is represented in Kunen-Tall \cite{KT79}, see history in the end of
\cite{Sh:80} and see more in \cite{Sh:546}.  We thank Mirna D\v{z}amonja for
drawing our attention to the problem and Ashutosh Kumar  
and Shimoni Garti 
for various
corrections and the referee for helpful suggestions.
\bigskip

\subsection {Are Some Versions of Axioms Equivalent?}\
\bigskip

To phrase our problem see the Definition below.
\newline
Kolesnikov asked:
\begin{question}
\label{x0}
Is there a forcing notion $\bbP$ satisfying
$(1)_a,(2)_b,(3)_{b,\omega}$ but not equivalent
to a forcing notion $\bbP'$ satisfying $(1)_a,(2)_b,(3)_a$?
\end{question}

\begin{definition}
\label{x2}
Consider the following conditions on a forcing notion $\bbP$ for a
   fixed  
   $\mu = \mu^{< \mu}$:
\medskip

\noindent
\underline{completeness}:

\noindent
$(1)_a \quad$ increasing chains of length $< \mu$ have a lub.
\smallskip

\noindent
$(1)_{a,< \theta} =
(1)_{a,\theta} \quad$ increasing chains of length $< \theta$ have a lub.
\smallskip

\noindent
$(1)_{a,\le  \theta} =
(1)_{a,\theta} \quad$ increasing chains of length $\le  \theta$ have a lub.
\smallskip

\noindent
$(1)_{a,=\theta} \quad$ increasing chains of length $\theta$ have a lub.
\smallskip

\noindent
$(1)_b \quad$ increasing chains of length $< \mu$ have a ub.
\smallskip

\noindent
$(1)_{b,< \theta} = (1)_{b,\theta} \quad$ increasing chains of
length $< \theta$ have an ub.
\smallskip

\noindent
$(1)_{b,\le  \theta} = (1)_{b,\theta} \quad$ increasing chains of
length $\le  \theta$ have an ub.
\smallskip

\noindent
$(1)_{b,=\theta} \quad$ increasing chains of length $\theta$ 
have an ub.  
\smallskip

\noindent
$(1)_c \quad \bbP$ is strategically $\alpha$-complete for every
$\alpha < \mu$, see \ref{z3}.
\smallskip

\noindent
$(1)_{c,\alpha} \quad \bbP$ is strategically $\alpha$-complete; 
where here
$\alpha \le \mu$.
\smallskip

\noindent
$(1)^+_c \quad$ there is a ``stronger" order $<_{\st}$ on $\bbP$ which means:
\mn
\begin{enumerate} 
\item[$\bullet_1$]  $p_1 <_{\st} p_2 \Rightarrow p_1 <_{\bbP} p_2$
\sn
\item[$\bullet_2$]  $p_1 \le_{\bbP} p_2 <_{\st} p_3 \le_{\bbP} p_4
  \Rightarrow p_1 <_{\st} p_4$
\sn
\item[$\bullet_3$]  any $<_{\st}$-increasing chain of length $< \mu$ has
 a $\le_{\bbP}$-ub (hence a $<_{\st}$-ub)
 \sn 
 \item[$\bullet_4$]  for every $ p $ there is $ q $ 
 satisfying $ p <_{\stt}  q $
\end{enumerate}
\smallskip

\noindent
$(1)_{d,< \theta} = (1)_{d,\theta} \quad$ any increasing continuous chain
of length $< \theta$ has a lub.
\smallskip

\noindent
$(1)_{d,=\theta}$ \quad any increasing continuous chain of length
$\theta$ has a lub.
\medskip

\noindent
\underline{Strong $\mu^+$-c.c.}: for a stationary $S \subseteq
S^{\mu^+}_\mu$, the default value being $S^{\mu^+}_\mu$, see \ref{z1};
we may write $(2)_x[S]$ when $S$ 
is 
is neither  
the default
value n or   
clear from the context.
\smallskip

\noindent 
\item[$(2)_a$]  
Given a sequence $\langle p_i:i < \mu^+\rangle$ of
 members of $\bbP$ there are a club $C$ of $\mu^+$ and a
regressive function $h$ on $C \cap S$
such that $\alpha,\beta \in
C \cap S \wedge h(\alpha) = h(\beta) \Rightarrow p_\alpha,p_\beta$
have a lub.

\noindent
\item[$(2)_b$] 
like $(2)_a$ but demanding just that $p_\alpha,p_\beta$
have an ub.
\smallskip

\noindent 
\item[$(2)^+_{a,\theta}$]  if $p_\alpha \in \bbP$ for $\alpha < \mu^+$
\then \, we can find a club $E$ of $\mu^+$ and a regressive $\hh:S \cap E
  \rightarrow \mu^+$ such that: if $i(*) < 1 + \theta,\delta_i \in S \cap
E$ for $i<i(*)$ and $\hh \rest \{\delta_i:i < i(*)\}$ is constant then
  $\{p_{\delta_i}:i < i(*)\}$ has a lub

\noindent 
\item[$(2)^{+}_{b,\theta}$]  like $(2)^+_{a,\theta}$
but in the end the set has a ub  

\noindent 
\item[$(2)^*_{a,\theta}$]  if $p_\alpha \in \bbP$ for $\alpha < \mu^+$
\then \, we can find $\bar q,E,\hh$ such that
\sn
\begin{enumerate}
\item[$\bullet_1$]  $\bar q = \langle q_\alpha:\alpha < \mu^+   
\rangle$
\sn
\item[$\bullet_2$]  $p_\alpha \le_{\bbP} q_\alpha$
\sn
\item[$\bullet_3$]  $E$ a club of $\mu^+$
\sn
\item[$\bullet_4$]  $h$ is a regressive function on $S \cap E$
\sn
\item[$\bullet_5$]  if $\cU \subseteq S \cap E$ 
has cardinality $<
 1 +
 \theta$ and $\hh \rest \cU$ is constant,
then $\{q_\delta:\delta \in \cU\}$ has a lub.
\end{enumerate} 

\noindent 
\item[$(2)^{*}_{b,\theta}$]  like $(2)^*_{a,\theta}$
but in the end the set has a ub  
\mn
\medskip

\noindent
\underline{For $\varepsilon < \mu$ a limit ordinal, .e.g. $\omega$}:
\smallskip

\noindent
\item[$(3)_a$]   
any two compatible $p_1,p_2 \in \bbP$ have a lub.
\smallskip

\noindent
\item[$(3)_{b,\varepsilon}$]  
if $\langle p_{\ell,\zeta}:\zeta <
\varepsilon\rangle$ is increasing for $\ell=1,2$ and
$p_{1,\zeta},p_{2,\zeta}$ are compatible for every $\zeta <
\varepsilon$ \then \, $\{p_{\ell,\zeta}:\ell
\in \{1,2\},\zeta < \varepsilon\}$ has an upper bound; recall
$\boxplus$ of \S(0A).

\noindent 
\item[$(3)_{b, \theta , \varepsilon }$] 
if  
(a) then (b) where:   
\begin{enumerate}
\item[(a)] 
\begin{enumerate} 
\item[$\bullet _1$]
$ p_{\zeta, i } \in \mathbb{P} $  for $ \zeta < \varepsilon $ and 
$ i < < i_* < \theta $ 
\item[$\bullet _2$] 
 if $ i < i_* $  then the sequence 
$\langle p_{\zeta, i}: \zeta < \varepsilon \rangle $ 
is $ < _{\stt}$-increasing ; (usually $< _{\stt } $ is from 
$ (1)^+_c $ 
\item[$\bullet _3$]
for each $ \zeta < \varepsilon $ the set
$ \{p_{\zeta, i }: i < i_* \}  $ has a common upper bound
\end{enumerate} 
\item[(b)]  the set 
$\{p_{\zeta, i }: \zeta < \varepsilon , i < i_* \} $
has a common upper bound.

\noindent 
\item[$(3)_{a, \theta, \varepsilon }$] 
 like $(3)_{b, \theta, \varepsilon }$ but in $ \bullet _3$  
 we have lub. 
\end{enumerate} 
\end{definition}


\begin{definition}
\label{x3}
Assume first  $D$ a normal filter on $\mu^+$ to which $S^{\mu^+}_\mu$ belongs
(we may omit $D$ when it is (the club filter on $\mu^+$) + $S^{\mu^+}_\mu$, see
Definition \ref{z6} we may omit $D$ if clear from the context).
We may write $ S $ instead $ D $  
when $ D $  is (the club filter on $\mu^+$) +
$S^{\mu^+}_\mu$.
Second 
$ 2 \le \theta \le \mu $, we may omit $ \theta $ when $ \theta = 2 $; 
we may write $ =\theta \text{ or }  \le \theta $  
instead $ \theta ^+$. or (essentially equivalent) 
$ \theta +1$
Third assume    
$\bbP$ is a forcing notion and  $\varepsilon < \mu$ is an ordinal;
a limit
ordinal if not said otherwise. 
Writing $ \xi  $ instead $ \varepsilon $  means
{`}{ }for every limit ordinal $ < \xi $ .
Note that  
$ (2)^ \varepsilon _{c,D }$  is equal to 
$ *^ \varepsilon _ {\mu , D }$ of \cite{Sh:546}.

Then we define the following conditions on $\bbP$:

\mn
\begin{enumerate}
\item[$(2)^\varepsilon_{c, \theta ,D}= (2)_{c,\theta, D, \varepsilon }$] 
in the following game the COM player
  has a winning strategy: 
\sn
\begin{enumerate}
\item[$(a)$]  a play last $\varepsilon$-moves
\sn
\item[$(b)$]  in the $\zeta$-th move 
a triple  
$(\bar p_\zeta,\hh_\zeta , S_ \zeta )$ is chosen  
  such that:
\sn
\item[${{}}$]  $(\alpha) \quad \bar p_\zeta
= \langle p_{\zeta,\alpha}:\alpha \in S_\zeta\rangle$
\sn
\item[${{}}$]  $(\beta) \quad p_{\zeta,\alpha} \in \bbP$
\sn
\item[${{}}$]  $(\gamma) \quad S_\zeta \in D$
\sn
\item[${{}}$]  $(\delta) \quad S_\zeta \subseteq \cap\{S_\xi:\xi < \zeta\}$
\sn
\item[${{}}$]  $(\varepsilon) \quad$ if $\alpha \in S_\zeta$ then $\langle
  p_{\xi,\alpha}:\xi \le \zeta\rangle$ is a $\le_{\bbP}$-increasing sequence
\sn
\item[${{}}$]  $(\zeta) \quad \hh_\zeta$ is a pressing down function on
  $S_\zeta$  
\sn
\item[$(c)$]  $\COM$ chooses\footnote{Why $1 + \zeta$ not, e.g. $\zeta
    +1$?  First, we like the INC to have the first move so that if
 $\bbP$ satisfies the condition and $p \in \bbP$ then $\bbP \rest
 \{q:p\le_{\bbP} q\}$ satisfies the condition.  Second, we like the
 player COM to move in limit stages, as this is a weaker demand.}
 $(\bar p_\zeta,h_\zeta)$ when $1 + \zeta$
  is even, $\INC$ chooses it when $1 + \zeta$ is odd
\sn
\item[$(d)$]  $\COM$ wins a play when it always could have
made a legal move, and in the end 
there is $S_\varepsilon \in D$
included in $\bigcap\limits_{\zeta < \varepsilon} S_\zeta$ such that:

if $i_* < \theta $ and $\alpha_i \in S_\varepsilon$
for $ i < i_*$
 and 
 for each $ i < i_* $  we have 
 $\bigwedge\limits_{\zeta < \varepsilon}
  \hh_\zeta(\alpha_i  
  ) = \hh_\zeta(\alpha_0)$ then 
  the set  
  $\{p_{\alpha_i ,\zeta}:\zeta < \varepsilon,
  i < i_* $ has 
  an ub
\end{enumerate}
\sn
\item[$(2)^\varepsilon_{d, \theta ,D}$] 
is defined as above replacing clause
 $(b)(\varepsilon)$ by
\sn
\item[${{}}$]  $(\varepsilon)' \quad$ if $\alpha \in S_\zeta$ then
  $\langle p_{\xi,\alpha}:\xi \le 
  \zeta\rangle$ is
  $\le_{\bbP}$-increasing continuous.
\end{enumerate}
\end{definition}

\begin{remark}
\label{x4}
1) So for a forcing notion $\bbQ,(2)^\varepsilon_{c,D}$
for $\varepsilon$ limit is 
$*^\varepsilon_D [\mathbb{Q} ] $ 
\cite[7]{Sh:546}.
.  
Also $ \mathbb{Q} $ satisfies 
$ (1)_b +(2)^2_{b,2,D} + (3)_a $
means  
$*^1_{\mu, \mathbb{Q} }$  
from the beginning of \S0(A).
Also $ \mathbb{Q} $  satisfies 
$ (1)_c +(2)^1_{a,2}$   mean   $ *^2_{\mu, \mathbb{Q} }$.

\noindent
2) Note that ``$\bbP$ satisfies
$(2)^\varepsilon_{c,D}$" implies
a weak version of strategic completeness (see  
  $ (1)_{b, \theta } $ for $ \theta = | \varepsilon |^+$).  
\end{remark}

\begin{definition}
\label{x5}
1) For suitable $ x,y,x $,  
  (but we may omit e.g. $ (3)_z$)
let $\Ax_{\lambda,\mu}((1)_x,(2)_y,(3)_z)$ 
mean: if $\bbP$ is a
   forcing notion satisfying those conditions and $\cI_i \subseteq
   \bbP$ is dense open for $i < i(*) < \lambda$ \then \, some directed
   $\mathbf G \subseteq \bbP$ meets every $\cI_i$.

\noindent
2) We may omit $\lambda$ if $\lambda = 2^\mu
\le \mu ^+$, we may  
more generally 
write
$\Ax_{\lambda,\mu}(K)$ for $K$ a property of forcing notion.

\noindent
3) For an ordinal\footnote{really omitting $(1)_b$ does not make a real
  difference but is natural} $\varepsilon < \mu$, a limit ordinal if
not said otherwise, let
$\Ax^\varepsilon_{\lambda,\mu}$ means:
$\Ax_{\lambda,\mu}( 
   (1)_c + (2)^\varepsilon_c)$,
we may omit
$\lambda$ if $\lambda = 2^\mu > \lambda^+$.
\end{definition}

\noindent
See on more axioms Roslanowski-Shelah \cite{Sh:655} parallel to
forcing and \cite{Sh:589} and references there.
In \S1 if we replace $C_\delta$ by a
stationary, co-stationary subset of $\delta$, we can iterate
appropriate $\mu^+$-c.c. $(< \mu)$-complete forcing notion. 
Earlier  
we have wondered (
for answers on this question see \ref{x11}(2)
\begin{question}
\label{x8}  Assume $ \mu  = \mu ^{< \mu }$
 
 1) In \cite{Sh:80},  can  the demand {`}{`}well met" cannot   be omitted?

\noindent 
2) Is there an example $\bbP$ where 
$ (1)_c +  (2)^\theta _c$  
 holds but
$  (1)_c +  
 (2)^ \partial _c$  
fails for any $\partial \in 
\Reg   
\backslash
\{\theta\}$ where $\theta = \cf(\theta),
\cf(\partial )   
= \partial <
\mu$?  The case $\partial = \aleph_0 < \theta$ is natural.  

\noindent
3) Do we have an example for $\Ax((1)_b + (2)_b + (3)_a)$ but not
 $\Ax^\varepsilon_\mu$  with e.g.   
 $ \varepsilon = \omega$, ?

\end{question}

\begin{discussion}  \label{x11}
1) Note: if we have $(3)_a =$ called well met then we have $(2)_a \equiv
(2)_b$.  If in addition we have $(1)_b$ \then \, we have
$(2)^\varepsilon_c$ for every $\varepsilon$.
Hence \ref{x8}(2) may be the true question.

\noindent 
2)
In \S1 (see \ref{a20}) we shall show that  the demand 
{`}{`}well met" cannot be omitted in \cite{Sh:80}; in other words,
the statement $ \Ax_\mu ((1)_a, (2)_b)$  is inconsistent.

In \S2 
for $ \theta ,  \partial <  \mu $  regular not equal 
we get  $ \Ax_\mu ((1)_c + (2)_{a, = \theta }$ 
   but 
   not $ \Ax _ \mu  ((1)_c + (2)_{a, \partial } ) $  
see \ref{c28}, but this does not answer 
Question \ref{x8}(2).
In \S3 we answer \ref{x8}(2).

\noindent
3) Suppose we consider a forcing notion as in \S1, i.e. for \S2 use
   $\theta=1$, but as in \ref{d6}, for $\alpha \in C_\delta \cap
   S^{\mu^+}_\theta$ no uniformization 
   is demanded.  
   This makes
   $\Ax^\theta_\mu$ holds for this forcing notion, but $*^\partial_\mu$
   fail, so all seems fine.

\noindent 
4) Below, in  
fact for $\langle C_\delta,\mathbf f_\delta:\delta \in S\rangle$,
   we may force also the $C_\delta$ (in $\bbQ$ in \S1); we may not ask
   that $C_\delta$ is closed in $\delta$ and let $\bar\alpha^*_\delta = \langle
   \alpha^*_{\delta,\xi}:\xi < \mu\rangle$ list $C_\delta$ in
increasing order so with limit $\delta$, but generically we can have
   $\alpha^*_{\delta_1,\zeta} = \alpha^*_{\delta_2,\zeta},\mathbf f_{\delta_1}
(\alpha^*_{\delta_1,\zeta}) \ne \mathbf
f_{\delta_2}(\alpha^*_{\delta_2,\zeta})$ for $*^1_\mu$, i.e. anyhow
   seems reasonable.
\end{discussion}

\begin{observation}
\label{x20}  
Assume $\mu = \mu^{< \mu}$ and $\varepsilon < \mu$ limit.

\noindent
1)For every limit
$\varepsilon < \mu$,  if  
the forcing notion $\bbQ$ satisfies the conditions
$(1)_{b,|\varepsilon|^+},(3)_a$ and $(2)_b$, here equivalently $(2)_a$
\then \, $\bbQ$ satisfies $*^\varepsilon_\mu$ of \cite{Sh:546},
i.e. $(2)^\varepsilon_c$ from Definition \ref{x3}, .

\noindent
2) If $\bbP$ satisfies $(3)_a$ \ then \, $\bbP$ satisfies
$(3)_{a,\varepsilon}$.

\noindent
3) If $\bbP$ satisfies $(1)_{b,|\varepsilon|^+} + (2)^+_{a,2}$ \then
\, $\bbP$ satisfies $(2)^\varepsilon_c$.

\noindent
4) For any $\bbP$ we have: $(1)_a \Rightarrow (1)_b \Rightarrow
(1)^+_c \Rightarrow (1)_c$ and $(1)_a \Rightarrow (1)_{d, \mu}  
\Rightarrow
(1)_c$. 
Similarly $(1)_{a,\theta}
\Rightarrow (1)_{b,\theta} \Rightarrow 
(1)_{c,\theta}$ and
$(1)_{a,=\theta} \Rightarrow (1)_{b,=\theta}$ and $(1)_{a,\theta}
\Rightarrow (1)_{d,\theta}$ and $(1)_{a,=\theta} \Rightarrow
(1)_{d,=\theta}$.

\noindent
5) For any $\bbP$ we have $(2)^{+}  
_{a,\theta} \Rightarrow
(2)^*_{a,\theta} \Rightarrow (2)^+_{b,\theta}$.

\noindent
6) If $\bbP$ satisfies $(2)^\varepsilon_{c,D}$ \then \, forcing with
$\mathbb{Q} $ adds  
no new sequence of ordinals of length $\le \varepsilon$.
\end{observation}

\begin{PROOF}{\ref{x20}}
Just read the definitions carefully.

E.g.

\noindent
3) Recall $\boxplus$ of \S(0A).
\end{PROOF}

\begin{claim}
\label{x24}
1) $\Ax^\varepsilon_\mu$, i.e. $\Ax_\mu(  
(2)^\varepsilon_c)$ is
equivalent to the axiom in \cite{Sh:546}.

\noindent
2) $\Ax_\mu((1)_b,(2)_a,(3)_a)$ is the axiom from \cite{Sh:80}.
If $\theta,\sigma$ are regular cardinals $< \mu$ and
$\Ax^\theta_\mu$ does not imply $\Ax^\sigma_\mu$ \then \,
$\Ax_\mu((1)_b,(2)_a,(3)_a)$ so the axiom from \cite{Sh:80}, does not imply
$\Ax^\sigma_\mu$.
\end{claim}

\begin{PROOF}{\ref{x24}}
Easy, too.
\end{PROOF}

\noindent
Many works on forcing for uniformizing see
\cite{Sh:64}, \cite{Sh:587}, \cite[Ch.VIII]{Sh:f} and on ZFC results
see \cite{Sh:65}, \cite[AP,\S1]{Sh:f}.
\bigskip

\subsection {Preliminaries} \
\bigskip

\begin{notation}
\label{z1}
For  regular $\theta < \lambda$ let $S^\lambda_\theta = \{\delta <
\lambda:\delta$ has cofinality $\theta\}$.
\end{notation}

\begin{definition}
\label{z3}
1) We say that a forcing notion
$\mathbb P$ is strategically $\alpha$-complete \underline{when}
for each $p \in \mathbb
P$ in the following game $\Game_\alpha(p,\mathbb P)$ between the players COM
and INC, the player COM has a winning strategy.

A play lasts $\alpha$ moves; in the $\beta$-th move, first the player
COM chooses $p_\beta \in \mathbb P$ such that $p \le_{\mathbb P} p_\beta$
and $\gamma < \beta \Rightarrow q_\gamma \le_{\mathbb P} p_\beta$ and
second the player INC chooses $q_\beta \in \mathbb P$ such that $p_\beta
\le_{\mathbb P} q_\beta$.

The player COM wins a play if he has a legal move for every $\beta <
\alpha$.

\noindent
2) We say that a forcing notion $\mathbb P$ is $(< \lambda)$-strategically
complete \underline{when} it is $\alpha$-strategically complete for every
$\alpha < \lambda$.
\end{definition}

\begin{definition}
\label{z6}
For a filter $D$ on a set $I$
\mn
\begin{enumerate}
\item[$(a)$]  $D^+ = \{A \subseteq I:I \backslash A \notin D\}$
\sn
\item[$(b)$]  for $S \in D^+$ let $D+S = \{A \subseteq I:A \cup (I
  \backslash S) \in D\}$.
\end{enumerate}
\end{definition}

\begin{theorem}
 \label{z9}
 Assume $ \mu = \mu ^{< \mu }$ and $ D $ is a normal 
 filter on $ \mu ^+$ to which $ S^{\mu ^+}_\mu $ belongs;
  not that in $ \mathbf{V} ^ \mathbb{P} $ we interpret  
  $ D $ as the normal filter on $ \mu ^+$ it generate.
  Assume further that $ 2 \le \theta \le \mu $.
  \Then
\,    each of the following properties listed in (B)
    of forcing notions is preserved
    by $ ( < \mu )$-support iteration which mean clause (A)
    is satisfied; where:
   
   \begin{enumerate} 
   \item[(A)]  if $ \mathbf{q} = 
  \langle \mathbb{P} _ \alpha , \name{\mathbb{Q}}_ \beta :
  \alpha \le \lg(\mathbf{q}), \beta < \lg (\mathbf{q} )\rangle $
  is a $ (< \mu )$-support iteration.
  and for each $ \beta  < \lg(\mathbf{q})$ we have
  $ \Vdash _{\mathbb{P} _ \beta  } \text{`}\text{`} ( \name{ \mathbb{Q} }_ \beta $ 
  satisfies the property $ \Pr$"  
  \then \, the forcing notion
  $ \mathbb{P} _ \mathbf{q} = \mathbb{P} _{\lg( \mathbf{q} }$
  satisfies the property $ \Pr $.  
  \item[(B)] 
  the property  
  $ \Pr $ of forcing notion $ \mathbb{Q} $ 
  is one of the following:  
  \begin{enumerate}  
 \item[(a)]    
 the property   
 $  (2)^ \varepsilon _ {c, D}$   for some 
 limit ordinal  $ \varepsilon < \mu $
 \item[(b)] the property $ (1)_{c, \theta } $
 \item[(c)] the property $ (1)^+_{c, \theta } $
 \item[(d)] for each  limit ordinal $ \varepsilon < \mu $,
 the property 
 $ (2)^\varepsilon _{c, \theta, D} $ 
 \item[(e)] for each  limit ordinal $ \varepsilon < \mu $,
 the property 
 $ (2)^\varepsilon _{d, \theta, D} $
  \end{enumerate} 
  \end{enumerate} 
 \end{theorem} 
 
 \noindent 
 
 \begin{PROOF} {\ref{z9}}
Cases (b),(c) are well known. 

 CASE (a)
 
 This holds by \cite{Sh:546}
 
 CASE (d)
 
 See Shelah-Spinas \cite{Sh:1110}.  

 CASE (e)
 
 Similalry.
 
 \end{PROOF} 
 
 \newpage  
 
\section {On $\mu^+$-regressive-c.c.; an example}


First, we shall concentrate on the case $\mu$ is not
strongly inaccessible.
\begin{hypothesis}
\label{a2}
1) $\mu = \mu^{< \mu} > \aleph_0$.

\noindent
2) $S \subseteq  S^{\mu^+}_\mu = \{\delta < \mu^+:\cf(\delta) = \mu\}$
is stationary, the main case is $ S = S^{\mu ^+}_\mu $.
\end{hypothesis}

\begin{definition}
\label{a4}
$\bar C$ is an $S$-club system \when \, $\bar C = \langle C_\delta:\delta
\in S\rangle,C_\delta$ a club of $\delta$ of order type $\mu$.
\end{definition}

\begin{definition}
\label{a7}
1) We say $(\cW,\bar{\mathbf f})$ is an $(S,\bar C,\kappa)$-parameter or
just a  $(\bar C,\kappa)$-parameter \when \,:
\mn
\begin{enumerate}
\item[$(a)$]    $S \subseteq S^{\mu^+}_\mu$ is stationary; see \ref{a2}(2),
\sn
\item[$(b)$]  $\bar C$ is an $S$-club-system so we may omit $S$
\sn
\item[$(c)$]  $ \le \mu$ is $\ge 2$,
if $\kappa=2$ we may omit $\kappa$ and write $\bar C$
\sn
\item[$(d)$]  $\cW \subseteq \mu$; if $\cW = \mu$ we may omit $\cW$
\sn
\item[$(e)$]  $\bar{\mathbf f} = \langle \mathbf f_\delta:\delta \in S\rangle$
\sn
\item[$(f)$]  $\mathbf f_\delta:C_\delta \rightarrow \theta $ 
\end{enumerate}
\mn
2) For 
$(\cW,\bar{\mathbf f})$ an $(S,\bar C,\kappa)$-parameter we define a
forcing notion $\bbQ = \bbQ  
_{(\cW,\bar{\mathbf f},\bar C)}$ as follows:
\mn
\begin{enumerate}
\item[$(A)$]  $p \in \bbQ$ \Iff \, $p$ consists of
\sn
\begin{enumerate}
\item[$(a)$]  $v 
\in [S]^{< \mu}$
\sn
\item[$(b)$]  $h$  is  a function with domain $v$  
\sn
\item[$(c)$]  if $\delta \in v$ then $h(\delta)$ is a 
non-empty 
 bounded
  subset of $\mu$ closed in its supremum
\sn
\item[$(d)$]  if $\delta_1,\delta_2 \in v$ and $\alpha \in  
  C_{\delta_1} \cap C_{\delta_2}$ 
  and $\otp(\alpha \cap C_{\delta_\ell})
  \in h  
  (\delta_\ell)$ and $\otp(C_{\delta_\ell}\cap \alpha) \in \cW$ for
  $\ell=1,2$ \then \, $\mathbf f_{\delta_1}(\alpha) = \mathbf f_{\delta_2}(\alpha)$
\sn
\item[$(e)$]  if 
$\delta_1 \ne \delta_2 \in v$ and $\beta \in  
  C_{\delta_1} \cap C_{\delta_2}$
then   
for $ {\ell} = 1,2 $  there is $ \beta _{\ell}  \in h_p(\delta _{\ell} )$
satisfying  $ \otp( C_{\delta _{\ell} }\cap \beta ) \le  \beta _{\ell} $ 
\end{enumerate}
\sn
\item[$(B)$]  $p \le_{\bbQ} q$ \Iff \,:
\sn
\begin{enumerate}
\item[$(a)$]  $v_p \subseteq v_q$  
\sn
\item[$(b)$]  $\delta \in v_p \Rightarrow h_p(\delta) \trianglelefteq  
  h_q(\delta)$ .
\end{enumerate}
\end{enumerate}

\noindent 
3) if $ {\mathscr W } = \mu $ we may omit it.
\end{definition}

\begin{definition}
\label{a9}
Let $(\cW,\bar{\mathbf{f}}  )$ be a 
$(\bar C,\kappa)$-parameter and let 
$\bbQ =
\bbQ   
_{\cW,\bar{ \mathbf{f}}  ,\bar C}$.  

\noindent
1) For $p \in \bbQ$ let $g_p$ be the function  
\mn
\begin{enumerate}
\item[$(a)$]  with domain
\begin{equation*}
\begin{array}{clcr}
\{\alpha:&\text{ some } \delta \text{ witnesses } \alpha \in
\Dom(h_p) \text{ which means } \delta \in v_p,\alpha \in C_\delta, \\   
  &\otp(C_\delta \cap \alpha) \in
   h_p(\delta)$ and $\otp(C_\delta \cap \alpha) \in \cW\}  
\end{array}
\end{equation*}
\sn
\item[$(b)$]  for $\alpha \in \Dom(g_p)$ we have:  

\[ 
g_p(\alpha)   
= \big( h_p  
(\delta)\big) (\alpha) \text{ for every   
  witness } \delta \text{ for } \alpha \in \dom(h_p).  
\]
\end{enumerate}
\mn
2) Let $\name g$ be the $\bbQ$-name 
for  
$\cup\{g_p:p \in \name{\mathbf  
  G}\}$.

\noindent
3) Let $\name E_\delta = \name E_\delta[\bbQ]$ be the $\bbQ$-name 
for 
   $\cup\{h_p(\delta):p \in \name{ \mathbf G } ,\delta \in v_p\}$ and let  
$\name{\cW}_\delta = \{\alpha \in \name E_\delta:\otp(C_\delta \cap \alpha)
  \in \cW\}$.
\end{definition}

\begin{claim}
\label{a12}
Assume $(\cW,\bar{\mathbf f})$ is an $(S,\bar C,\kappa)$-parameter and
$\bbQ = \bbQ_{(\cW,\bar{\mathbf f},\bar C)}$.

\noindent
1) $\bbQ$ is $(< \mu)$-complete, moreover any $\le_{\bbQ}$-increasing
   sequence of length $< \mu$ has a $\le_{\bbQ}$-lub.

\noindent
2) If $\delta \in S$ and $\hi < \mu$ \then \, the following subsets
of $\bbQ$ are dense and open:
\mn
\begin{enumerate}
\item[$\bullet_1$]  $\cI_\delta = \{p \in \bbQ:\delta \in v_p\}$  
\sn
\item[$\bullet_2$]  $\cI_{\delta,\hi} = \{p \in \cI_\delta:\hi <
  \sup(h_p(\delta))\}$  
\sn
\item[$\bullet_3$]  $\cI^*_\hi = \{p \in 
\mathbb{Q}    
$: if $\delta \in v_p$ then $\hi  
  < \sup(h_p(\delta))
  \text{ and }h_p( \delta ) \text{ has a last member}
  \}$.   
\end{enumerate}
\mn
3) 
For every $ \delta \in S $, the function 
$\name g$ almost extends $\mathbf{f}  
_\delta$, i.e. $\Vdash_{\bbQ} ``\name g \supseteq  
\mathbf{f} _\delta \rest \{\alpha \in
C_\delta:\otp(\alpha \cap C_\delta) \in
\name \cW_\delta\}$, 
$ {\mathscr W } _ \delta = {\mathscr W } \cap \name{ E}_\delta $
and   $\name E_\delta$ 
is 
a club of $\mu$ 
and if $\cW =
\mu$ then $\name {\mathscr W }   
_\delta$ is a club of $\mu"$.
\end{claim}

\begin{PROOF}{\ref{a12}}
1) Straightforward, see clause (A)(e) of Definition \ref{a7}(2) in particular.

\noindent
2),3)  Also easy.
\end{PROOF}

\begin{claim}
\label{a14}
Let $(\cW,\bar{\mathbf f}),(S,\bar C,\kappa),\bbQ$ be as above.

\noindent 
\Then  \,
$\bbQ$ satisfies clause $(2)_b$ of Definition \ref{x2} that is:

\mn
\begin{enumerate}
\item[$*^0_\mu$]   if $\bar p = \langle
p_\alpha:\alpha \in S\rangle$ and $\alpha \in S \Rightarrow
p_\alpha \in \bbQ$ \then \, there is a club $E$ of
$\mu^+$ and pressing down function $f:S \cap E \rightarrow \mu^+$, i.e.
$f(\delta) < \delta$, such that:
$(\delta_1 \ne \delta_2 \in S \cap E) \wedge f(\delta_1) =
f(\delta_2) \Rightarrow p_{\delta_1},p_{\delta_2}$ are compatible.
\end{enumerate} 


\end{claim}

\begin{PROOF}{\ref{a14}}
First, by \ref{a12}(1)(2), we choose $\langle q_\alpha:
\alpha \in S\rangle$ such that, for every $ \alpha \in S$:  
\mn
\begin{enumerate}
\item[$\odot_1$]  $(a) \quad p_\alpha \le q_\alpha$
\sn
\item[${{}}$]  $(b) \quad$ if $\delta \in v_{q_\alpha}$ but  
$\delta > \alpha$ then $\otp(C_\delta \cap \alpha) <
\sup(h_{q_\alpha}(\delta))$  
\sn
\item[${{}}$]  $(c) \quad 
\alpha \in v_{q_\alpha}$.  
\snn 
\item[{{}}]  $(d)$  
 $ h_p(\alpha )$  has a last element.  
\end{enumerate}
\mn

Second, choose a club $E$ of $\mu^+$ such that $\alpha \in S \cap E
\Rightarrow \sup(v_{q_\alpha}) < \min((E \backslash (\alpha +1))$.  

Third, choose  a
regressive 
function $\hh$ with domain $E \cap S$ such that:  
\mn
\begin{enumerate}
\item[$\odot_2$]  if $\delta(1) = \delta_1 < \delta_2 = \delta(2)$
are from $E \cap S$ and
  $\hh(\delta_1) = \hh (\delta_2)$ and 
  $\langle \alpha_{\ell,i}:i <
\otp(v_{q_{\delta(\ell)}})\rangle$ lists  
$v_{q_{\delta(\ell)}}$ in 
  increasing order for $\ell=1,2$ \then \, for some $j_*$:
\sn
\begin{enumerate}
\item[$(a)$]  $\otp(v_{q_{\delta(1)}}) = \otp(v_{q_{\delta(2)}})$ call  
  it $i(*)$
\sn
\item[$(b)$]  $j_* < i(*)$ and $\alpha_{1,j_*} =
  \delta_1,\alpha_{2,j_*} = \delta_2$
\sn
\item[$(c)$]  if $j < j_*$ then $\alpha_{1,j} = \alpha_{2,j}$
\sn
\item[$(d)$]  if $j > j_*$ but $j < i(*)$ then $C_{\alpha_{1,j}} \cap
  \delta_1 = C_{\alpha_{2,j}} \cap \delta_2$
\sn
\item[$(e)$]  $h_{q_{\delta(1)}}(\alpha_{1,i}) =  
h_{q_{\delta(2)}}(\alpha_{2,i})$ for $i < i(*)$  
\sn
\item[$(f)$]  if $\varepsilon \in h_{q_{\delta(1)}}(\delta_1)$ then the 
  $\varepsilon$-th member of $C_{\delta_1}$ is equal to the
  $\varepsilon$-th member of $C_{\delta_2}$.
\end{enumerate}
\end{enumerate}
\mn
Now it suffices to prove:
\mn
\begin{enumerate}
\item[$\odot_3$]  if $\delta_1 \ne \delta_2 \in S \cap E$ and
  $\hh(\delta_1) = \hh(\delta_2)$ 
  then $q_{\delta_1},q_{\delta_2}$ are
  compatible in $\bbQ$, 
\end{enumerate}
\mn
Why?  Define $q$ as follows:
\mn
\begin{enumerate}
\item[$\bullet_1$]   $v_q = v_{q_{\delta(1)}} \cup v_{q_{\delta(2)}}$  
\sn
\item[$\bullet_2$]  $h _ q   
(\delta) = h_{q_{\delta(\ell)}}(\delta)$ if $\ell \in  
  \{1,2\}$  and $\delta \in v_q  
  \backslash \{\delta_\ell\}$ 
\sn
\item[$\bullet_3$]  $h_q(\delta_\ell) =
    h_{q_{\delta(\ell)}}(\delta_\ell) 
\cup \{\beta _{\ell} \}$ where $\beta _{\ell}  < \mu,
\beta _{\ell}  > \max\{h_{q_{\delta(1)}}(\delta_1) \cup   
h_{q_{\delta(2)}}(\delta_2)\}$  and $\gamma >  
\sup\{\otp(\alpha \cap C_{\delta_\ell}):\ell \in \{1,2\}$ and
$\beta _{\ell}   
\in C_{\delta_1} \cap C_{\delta_2}\}$.
\end{enumerate}
\mn

First, 
$ g $ is well defined because in $ \bullet _2 $,  if $ h_q( \alpha ) $   
is defined in two ways, then  
$ \alpha < \delta _1 $ and they are equal because 
of $ \odot _2 $ 

Second, why $ q \in \mathbb{Q} $? We have to 
to check clauses (a)-(e) of Def \ref{a7}(2)(A).
Now clauses (a),(b),(c) are obvious. For clause (d),  
assume $ \gamma _1, \gamma _2 \in v_q $, and 
$ \alpha \in C_{\gamma_1} \cap C_{\gamma _2}$ 
and $ \otp (C_{\gamma _ {\ell} }\cap \gamma _ {\ell})
\in h_q ( \gamma _ {\ell} ) \cap {\mathscr W } $  
for $ {\ell} = 1,2$.
If $ \gamma _1, \gamma _2 \in v_{q_{\delta (1)}}$ 
then use $ q_{\delta (1)}\in \mathbb{Q} $,
and similarly if
 $ \gamma _1, \gamma _2 \in v_{q_{\delta (2)}}$ 
then use $ q_{\delta (2)}\in \mathbb{Q} $.
So \wilog \, $ \gamma _1 \in v_{q _{\delta (1)}} \setminus
v_{q _{\delta (2)}}$
and
$ \gamma _2 \in v_{q _{\delta (2)}} \setminus
v_{q _{\delta (1) }}$, so necessarily  
$ \gamma _1\ge  \delta (1), \gamma _2 
\ge \delta _2$  and
$ \alpha \in C_{\gamma _1} \cap C_{\gamma _2}
\subseteq \delta _1  
\cap \delta _2 $ (using the choice of $ \bar{ \mathbf{c} }$); 
using the notation of $ \odot _2$
let $ i( {\ell} )$ 
be such that  $ \gamma_{\ell} = \alpha _{{\ell}, i( {\ell} }$
so $ i( {\ell} ) \in [ j(*), i(*))$  for $ i({\ell} )
= 1,2$.
Now we get the result by applying 
clause (d) for 
$  q_{\delta (2)} \in \mathbb{Q} $ 
for $ \gamma_1, \gamma _2, \alpha _{2, i(1)}$  
$ \alpha _{2, i(1)}, \alpha _{2, i(2)}= \gamma _2$ recalling
$ \odot(d),(e)$, noting that in the case
$ ( \gamma _1, \gamma _2 )= ( \delta_1, \delta_2)$  necessarily 
$ \alpha \not= \beta _1, \beta_2 $ so 
$ \otp(C_{\delta (1)} \cap \alpha ) = 
  \otp( C_{\delta (2)} \cap \alpha  \in 
    h_{p _{\delta ( 1)} }(\alpha ) = 
     h_{p_{\delta (2)}}$

We are left with clause (e) which is proved similarly, 
recalling $ \bullet _3$  above.



It is easy to check that $q \in \bbQ$ and $q_{\delta_1} \le q,
q_{\delta_2} \le q$, so $\odot_3$ holds indeed

   
   
   


\end{PROOF}

\begin{theorem} 
\label{a17}
If (A) then (B) where
\begin{enumerate} 
\item[(A)] $ \mu , S, \bar{ C } ,\kappa, \theta   $  satisfy
\begin{enumerate} 
\item[(a)]  $ \mu = \mu ^{< \mu }  > {\aleph_0} $
\item[(b)] $ S = S^{\mu +}_\mu $
\item[(c)]  $ \bar{ C }= \langle C_ \delta : \delta \in S  \rangle $ 
is an $ S $-club system 
and for $ \delta \in S, \alpha < \delta $ (mainly $ \alpha \in C_\delta $);
we let  
$ \eta _\delta \in {}^{ \mu } \delta $  
list $ C_ \delta $  in increasing order
\item[(d)] $ F $ is a function from $ {\mathscr F} _{\mu  }$
to $ \kappa $  where  $ {\mathscr F} _ \mu = 
\{ f: f \text{ is a function from some } u
\in [\mu ^+ ]^{< \mu } 
\text{ to }\mu \} $; 
the default case is $ F(f)= f(\max(\dom(f))$  when well defined 
and zero otherwise.
\item[(e)] $ \bar{ a } = \langle a_{\delta , \hi } :
  \delta \in S, \hi < \mu \rangle $ 
where $ a_{\delta, \hi } \subseteq \eta _ \delta (\hi )+ 1$;
the default value of $ a_{\delta, \hi }$ is $ \{ \eta _\delta ( \hi ) \} $
\item[(f)] 
either $ \mu $ is a (strongly)  inaccessible  cardinal,
   and $ \theta < \kappa = \mu $ \underline({or}
$  \kappa = 2, \theta 
< \mu = 2 ^\theta $ 
\end{enumerate}
\item[(B)] we can find $ \bar{ \mathbf{c} }$  satisfying:
\begin{enumerate}
\item[(a)] $ \bar{\mathbf{c}} = \langle \mathbf{c} _ \delta : 
\delta \in S\rangle $
\item[(b)] $ \mathbf{c} _ \delta $ is a function from $ C_ \delta $  
to $ \kappa $ 
\item[(c)]  if $ f $ is a function from $ \mu ^+ $ 
   to $ \kappa $, 
then for stationarily many $ \delta \in S $, for stationarily 
many $ \varepsilon \in C_ \delta  $  we have:

$ \kappa = 2 \Rightarrow \mathbf{c} _\delta (\alpha  ) = 
F( f \rest a_{\delta, \hi }) $ 

and   

$ \kappa = \mu \Rightarrow  
\mathbf{c} _\delta (\alpha ) \not= 
F( f \rest a_{\delta, \hi }) $  
\end{enumerate} 
\end{enumerate} 
\end{theorem} 

\begin{discussion}  \label{a18}
  See \cite[AP.3.9,pg.990]{Sh:f}. But
  there, only the case
$ \mu = {\aleph_1}, \kappa = 2 $ is really proved, 
the case $ \mu $ an accessible cardinal and $ \kappa = 2 $ 
is stated to be similar.
In the case $ \mu $ inaccessibe, $ \kappa =2$, the statement
consistently fail as said in \cite[3.8(1)]{Sh:f},
see \cite{Sh:64}, \cite{Sh:186} and \cite{Sh:587}.
So by a request
we give here a full proof. 
\end{discussion} 

\begin{PROOF}{\ref{a17}}
Why?


Let $\lambda$ be big enough (e.g., ($2^{{\mu ^+}})^+)$,
and $M^*$ be an expansion of
 $({\mathscr H} (\lambda),\in)$ by Skolem functions
(so countably many; essentially, if we expand just  by 
a definable well ordering it suffices).

Suppose   
toward contradiction that clause (A) hold but clause (B) fail. 
It is known that there is a function $G$ from
$\{A:A\subseteq {\mu ^+}, |A|  < \mu 
\}$ to $\mu $ such that
$G(A)=G(B)$ implies $A,B$ have the same order
type and their intersection is an initial segment of both
(e.g. if $h_\alpha:\alpha\rightarrow \mu $ is one-to-one for
$\alpha<\mu $, we let
$G_0(A)\eqdf\{(\otp (A\cap \alpha), \otp (A\cap \beta),
 h_\beta(\alpha)):
\alpha\in A\hbox{ and } \beta\in A\}$. Now $G_0$ is as required except that
$\Rang (G_0)\not\subseteq \mu $ but $|\Rang (G_0)|\le \mu $
so we can correct thi by renaming).  

We shall  
now define 
 for any $\hp\in {\mathscr H} (\lambda)$ a sequence  
$\langle  \hc^\hp_\delta:\delta\in S  \rangle$ where 
$\hc^\hp_\delta:\mu \rightarrow {\mathscr H} (\mu )$, which we shall use later.

For every $\delta\in S  $, $i<\mu $,
let $N^\hp_{\delta,i}$ be the minimal submodel of $M^*$
(so closed under Skolem functions)  
including $\{\delta,i,\hp\}$  
such that its intersection with $ \mu $ is an ordinal 
 and let
\begin{enumerate} 
\item[$(*)_1$] let
\begin{enumerate}
\item[(a)]  $ \pi ^{\hp} _{\delta, \hi}$ be the Mostowski
collapse mapping from $ N^\hp_{\delta, \hi}$
\item[(b)] $ \mathbf{c} ^\hp_{\delta }$ is a function 
from $ \mu $ into $ {\mathscr H} (\mu ) $
\item[(c)]  for $ \hi < \mu $ we let 
$\hc^\hp_\delta(\hi)\eqdf\langle $ 
($(^\hp_{\delta , \hi }(N^\hp_{\delta,\hi},\hp,\delta,\hi),
G (N^\hp_{\delta,\hi}\cap {\mu ^+})\rangle$
  which belongs to $ {\mathscr H} ( \mu )$.
\end{enumerate} 
\end{enumerate} 

Note that 
\noindent
$(N^\hp_{\delta,i},\hp ,i,\delta)$ is $N^\hp_{\delta,i}$
expanded by three individual constants.

Now recall that toward contradiction   we are  assuming 
that clause (B) of the theorem fail.
This means that 

\begin{enumerate} 
\item[$(*)_2$] 
for every  sequence 
$\bar{ \hc }= \langle \hc_\delta:\delta\in S  \rangle $ 
where $ \hc_ \delta $ is a function from $ C_ \delta $ to $ \kappa $
there is
$h_{\hc}   
:{\mu ^+}\rightarrow \kappa $ 
such that:

for a closed unbounded
 set
of $\delta\in S  $, for a closed unbounded set of 
$ \hi < \mu $ we have  
$\hc_\delta (\hi  ) 
=F(h_{\bar{ \hc}}\uhr a_{\delta, \hi }$;
note that in the case $ \kappa = 2 $, replacing non-equal by equal make no difference!
\end{enumerate} 

Now  

\begin{enumerate} 
\item[$ (*)_3$]  in $ (*)_2 $ 
    we   can replace $ \kappa $ by the set 
$ {\mathscr H} ( \mu ) $, by changing $ F $  
\end{enumerate} 
[Why? If $ \kappa = \mu $  this is obvious as 
$ \mu $ and $ {\mathscr H} ( \mu ) $  have the same cardinality. 
So we can assume $ \kappa = 2$, and we can replace $ {\mathscr H} ( \mu )$  
by $ ^\theta 2$ because the later has cardinality $ \mu $.
For  $ \hn < \theta $ and 
$h$ any function into $^ \theta 2$, let $h^{[\hn]}$ be
defined by 
$h^{[\hn]}(\alpha )=(h(\alpha ))(\hn)$ for $\alpha \in \Dom (h)$.
Define the function $F^*$ by: $F^*(h)=\langle  F(h^{[\hn]}):\hn<\theta \rangle$.
We shall prove that replacing $ F $ by $ F^* $, the statement
$ (*)_3$  holds.
So assume 
we are given $\langle  \hc_\delta: \delta\in S  \rangle$ where  
$\hc_\delta\in {}^{\mu }({}^\theta 2)$, i.e.,
$\hc_\delta:  \mu   
\rightarrow {}^\theta 2 $; 
then for $ \hn < \theta $
the function $\hc^{[\hn]}_\delta\in {}^{\mu }2$ is well defined for each
$\delta\in S  $. 
Now for each $ \hn < \theta $, we can apply $ (*)_2$ 
so we can choose $h^{(\hn)}:{\mu ^+}\rightarrow 2$
such that for a club of $\delta\in S  $ for a club of
$ \hi < \mu $ 
we have 
$$\hc^{[\hn]}(\hi )=
   F(h^{(\hn)}(\hi )\uhr a_{\delta, \hi })  $$   
\hz
Define $h:{\mu ^+} \rightarrow {}^\theta  2$ 
by $h(\alpha )=\langle 
 h^{(\hn)}(\alpha ):\hn<\theta \rangle$, it is as required. 
So $ (*)_3 $ holds indeed.]
\hz



Now we
shall 
define by induction on $\hn<\theta $, $\hp(\hn)\in {\mathscr H} (\lambda)$, and
$h_\hn:{\mu ^+}\rightarrow {\mathscr H} (\mu )$.

Arriving to $ \hn $, let 
$ \hp(\hn ) = ( \langle h_\zeta : \zeta < \hn \rangle , \bar{C}, F, \bar{ a } )$
and let $\hc^{\hp(\hn)}_\delta:\mu \rightarrow {\mathscr H} (\mu )$
be as we have defined before (in $(*)_1$), 
so by $(*)_3$ 
\begin{enumerate} 
\item[$(*)_4$] 
there   are $ h_ \hn , W_\hn, \bar{ W}_ \hn $  such that:  
 \begin{enumerate}   
\item[(a)]     $h_\hn:{\mu ^+}\rightarrow {\mathscr H} (\mu )$;
\item[(b)]  $W^\hn\subseteq {\mu ^+}$  is a closed unbounded 
    subset of $ \mu ^+$
\item[(c)] $ \bar{ W }_ \hn = \langle W^ \hn _ \delta :
  \delta \in W \cap S  \rangle $
 \item[(d)] for every $\delta\in W^\hn\cap S  $,
$W^\hn_\delta\ $ 
is a closed unbounded subset of $ \mu $ 
 \item[(e)] for
$ \hi \in W^\hn_\delta, \delta\in W^n\cap S  $ we have:
 $\hc^{\hp(\hn)}_{\delta}(\hi ) =  
 F^*(h_\hn\uhr a_{\delta, \hi }) $ 
\end{enumerate} 
\end{enumerate} 


Now 

\begin{enumerate} 
\item[$(*)_5 $]  let  
\begin{enumerate} 
\item[(a)]  let $W=\bigcap_{\hn< \theta }W^\hn$, 
\item[(b)]  for $\delta\in W$ let 
$W_\delta=\bigcap_{\hn<\theta }W^\hn_\delta$.
\end{enumerate} 
\end{enumerate} 
Clearly $W$ is a closed unbounded subset of ${\mu ^+}$, and 
$W_\delta$ is a closed unbounded subset of $ \mu $ 
  for $ \delta \in W \cap S  $.  
So for every $\delta\in W\cap S  $, we can choose 
$\hi(\delta)\in W_\delta;$
hence 
by Fodor lemma, for some
$\hi(*)<{\mu ^+}$ and $\alpha (*)<\mu^+$ 
and $ \nu, \bar{ b } $  the set
$S_* = 
\{\delta\in W\cap S  :  
\hi(\delta)=\hi(*), 
\alpha ( \delta ) = \alpha (*), \eta _ \delta \rest (\xi +1) = \nu , 
\langle a_{\delta, i}: i \le \hi (*) \rangle = \bar{ b}
\}$ is stationary. 
As $\mu = \mu ^{< \mu }$   holds
there are $\delta_1, \delta_2$ 
and $\xi <\mu  $
such that:
\begin{enumerate}
\item[$(*)_6$]
\begin{enumerate} 
\
\item[(A)] $\delta_1<\delta_2$ are from $ S_*$
\item[(B)] $\xi \in W_{\delta_\ell}$ for
         $\ell=1,2.$
\item[(C)] $\eta_{\delta_1}(\xi)=\eta_{\delta_2}(\xi)$ 
\item[(D)] $\eta_{\delta_1}\uhr (\xi+1)=\eta_{\delta_2}\uhr (\xi+1)$        
\item[(E)] $ \langle a_{\delta _1, \hi }: \hi \le \hi(*) \rangle =
   \langle a_{\delta _2, \hi }: \hi \le \hi(*) \rangle $ 
\end{enumerate} 
\end{enumerate} 

So clearly we can assume
\begin{enumerate} 
\item[$(*)_7$] there are no $\delta^\dagger_1, \delta^\dagger_2$
         satisfying (A)-(E) 
         such that
         $\delta^\dagger_1\le \delta_1$, $\delta_2^\dagger\le \delta_2$
 and
         $(\delta_1^\dagger, \delta_2^\dagger)\neq (\delta_1,
 \delta_2)$.
\end{enumerate} 

Now as $\delta_1<\delta_2$,
for some $\hi>\xi, \eta_{\delta_1}(\hi)\neq \eta_{\delta_2}(\hi)$, and there
 is
a minimal such $\hi$;
but as $\eta_{\delta_1}, \eta_{\delta_2}$ are increasing and continuous,
such minimal $\hi$ should be a successor ordinal,
so 
similarly 
there is a maximal $\zeta$ among those satisfying
 $\zeta<\mu , \eta_{\delta_1}\uhr \zeta=\eta_{\delta_2}\uhr\zeta$,
$\eta_{\delta_1}(\zeta)=\eta_{\delta_2}(\zeta)$
and $\zeta\in W_{\delta_1}\cap W_{\delta_2}$. 


So

\begin{enumerate} 
\item[$(*)_8$]\quad $\hc^{\hp(\hn)}_{\delta_1} (\zeta)= \hc^{\hp(\hn)}_{\delta_2}(\zeta)$ for every $\hn  < \theta $
\end{enumerate} 

\nz
[Why?
as both are equal to $F^*(h_\hn\uhr a_{\delta_\ell ,
 \zeta})$.]  
 
 Fix a non-zero $ \hn < \theta $ for a while.  
 Looking at the definition of
 $\hc^{\hp(\hn)}_\delta(\zeta)$ (see $(*)_1 $) 
 we see that
$N^{\hp(\hn)}_{\delta_1, \zeta}$ is isomorphic to $N^{\hp(\hn)}_{\delta_2,
 \zeta}$,
and let the isomorphism be called $g_\hn$.
Note that the isomorphism is unique (as $\in$ in those models is
 transitive  and 
well founded 
and maps $ \bar{ C }, F, \bar{ a } $  to themselves.).

By the definition of $\hc^{\hp(\hn)}_\delta(\zeta)$,
clearly 
\begin{enumerate} 
\item[$(*)_9 $] 
\begin{enumerate} 
\item[(a)] $g_\hn(\hp(\hn))=\hp(\hn)$  hence 
     $ g_ \hn ( ( \bar{ C }  , F, \bar{ a } ))-=  
     ( \bar{ C }  , F, \bar{ a } ) $ 
\item[(b)] $  g_\hn(\delta_1)=\delta_2, g_\hn(\zeta)=\zeta, g_ \hn (\hn )=\hn$
\item[(c)] $g_\hn(\eta_{\delta_1})=\eta_{\delta_2}$
\item[(d)]  $g_\hn(W^{\hn-1})=W^{\hn-1}$ and
$g_\hn(W^{\hn-1}_{\delta_1})=
  W^{\hn-1}_{\delta_2}$ 
\item[(e)] $g_\hn(N^{\hp(\hn-1)}_{\delta_1, \zeta})=
   g_ \hn (N^{\hp(\hn-1)}_{\delta_2, \zeta})\in
N^{\hp(\hn)}_{\delta_2, \zeta}$.
\end{enumerate} 
\end{enumerate} 
[Why? Look at the definition of $ \hp ( \hn )$]

As $N^{\hp(\hn-1)}_{\delta_\ell, \zeta}$ is 
of cardinality $ < \mu $, its intersection with 
$ \mu $ is an ordinal 
and it belongs to
$N^{\hp(\hn)}_{\delta_\ell, \zeta}$, it is also included in it,
hence $g_\hn\uhr N^{\hp(\hn-1)}_{\delta_1, \zeta}$ is an isomorphism from
$N^{\hp(\hn-1)}_{\delta_1, \zeta}$ onto $N^{\hp(\hn-1)}_{\delta_2, \zeta}$ hence
(by the uniqueness of $g_\hn$  and $ (*)_9(b)$):
\begin{enumerate} 
\item[$(*)_{10}$]   
\quad $g_\hn\supseteq g_{\zeta }$  for $ \zeta < \hn $.
\end{enumerate} 

We now stop fixing $ \hn$.
For $\ell=1,2$, (recalling $ \theta < \mu $) let
$N_\ell=\bigcup_{\hn<\theta }N^{\hp(\hn)}_{\delta_\ell, \zeta}$ and 
$g=\bigcup_{\hn < \theta } g_\hn$;
so $g$ is an isomorphism from $N_1$ to $N_2$.
\nz
By the definition of $\hc^{\hp(\hn)}_{\delta_\ell}(\zeta)$, clearly the
second coordinates are the same, thus:
\begin{enumerate} 
\item[$(*)_{11}$]   
\quad
$G(N^{\hp(\hn)}_{\delta_1, \zeta}\cap {\mu ^+})
=G(N^{\hp(\hn)}_{\delta_2, \zeta}\cap{\mu ^+})$,
\end{enumerate} 

Hence those sets have their intersection an initial segment of both;
 as this holds for every $ \hn < \theta $,  clearly 
$N_1\cap {\mu ^+}, N_2\cap {\mu ^+}$ have their intersection
 an
initial segment of both
(as usually, we are not strictly distinguishing between a model and its
universe),
hence $g$ is the identity on $N_1\cap N_2\cap{\mu ^+}$.

Note that clearly $\delta_1\notin N_2$ as
 $g(\delta_1)=\delta_2\neq\delta_1$,
hence $\delta_2\notin N_1$. Now

\begin{enumerate} 
  \item[$ (*)_{12}$]
\item[(a)]  Letting  $\delta^*_\ell\eqdf \Min ({\mu ^+}\cap N_\ell\setminus(N_1\cap
 N_2))$, we have: 
$\delta^*_\ell\le \delta_\ell$, is a limit ordinal 
\item[(b)] $g(\delta_1^*)=\delta_2^*$ and so 
\item[(c)] $\cf (\delta_1^*)=\cf (\delta^*_2)$.
\item[(d)]   
   \quad $\cf (\delta^*_\ell)=\mu $.
\end{enumerate} 

Why? 
Clauses (a),(b) are obvious and use (c) follows. 
Clause (d) (that is $ \cf(\delta ^*_{\ell} )$)  
 holds as  
otherwise 
for some regular cardinal $ \sigma < \mu $ we have 
$\cf(\delta^*_1)=\sigma   $,
and as $\delta^*_1\in N_1$ for some $\hn < \theta $,
$\delta_1\in N^{\hp(\hn)}_{\delta_1, \zeta}$, hence there is
$\{\beta_\hm:\hm<\sigma  \}\in  
\delta^*_1\cap N^{\hp(\hn)}_{\delta_1,\zeta}$
cofinal in $\delta^*_1$. 
As $ \sigma < \mu $ necessarily it is included in 
$ N^{\hp(\hn)}_{\delta_1,\zeta} $, \wilog \,
$ \beta _ \iota $ is increasing with $\iota $. 
By the choice of $\delta^*_1$, if $ \iota < \sigma  $ then 
$ \beta_\iota \in N_1\cap N_2$,
hence $g(\beta_\hm)=\beta_\hm$;
let $\beta^*=\min(N^{\hp(\hn)}_{\delta_2,\zeta}\setminus \bigcup_\hm
 \beta_\hm)$,
so $\beta^*\in N^{\hp(\hn)}_{\delta_2, \zeta}\subseteq N^{\hp(\hn+1)}_{\delta_2,
 \zeta}$,
so $\delta^*_1
= \sup\{\beta_\iota :\iota <\sigma \}
=\sup(\beta^*\cap N^{\hp(\hn)}_{\delta_2, \zeta})\in N_2$,
contradiction.
\hz
So we have proved $(* )_{12}$.]
\hz

Now let for $\ell=1,2,\hi_\ell\eqdf N_\ell\cap \mu $,
(this intersection   
is an initial segment of $ \mu$) and $\beta_\ell\eqdf
\sup(N_\ell\cap \delta_\ell^*)$
hence $\beta_1=\beta_2$
(by $\delta^*_\ell$ definition) and call it $\beta$.
As $\cf(\delta^*_\ell)= \mu $ clearly
$\delta^*_\ell\ge \mu $, and so
clearly by $g$'s existence $\hi_1=\hi_2$ and call it 
$ \hi(*) $.  
(also as $\mu \in N_1\cap N_2\cap {\mu ^+}$, necessarily $N_1\cap
\mu  = N_2\cap \mu $).

As $\eta_{\delta^*_1}$ is a one to one function (being increasing)
from $\mu $,
clearly
\begin{enumerate} 
\item[$(*)_{13}$]   for every $ \hi < \mu $  we have 
$\eta_{\delta^*_1}(\hi)\in N_1 \iff \hi<\hi(*).$
\end{enumerate} 
Also $N_1\models \lq \langle  \eta_{\delta^*_1}(\hi):
\hi<\mu \rangle \rq$ is unbounded below
$\delta^*_1 $ (remember $N_1\prec M^*$ as
$N^{\hp(\hn)}_{\delta_1, \zeta}\prec M^*$ for each $\hn$).

So clearly $\beta=\sup\{\eta_{\delta_1^*}(\hi):\hi<\hi(*)\}$;
but $\eta_{\delta^*_1}$ is increasing continuous 
and $\hi(*)$ is a limit
ordinal (being $N_\ell\cap \mu $), hence
 $\beta=\eta_{\delta^*_1}(\hi)$.

 For the same reasons $\beta=\eta_{\delta_2^*}(\alpha)$.

Similarly   
 $\eta_{\delta^*_1}\uhr \alpha =\eta_{\delta^*_2}\uhr \alpha$
 because
 $g(\eta_{\delta^*_1})=\eta_{\delta^*_2}$,
 and $\alpha\in W^\hn_{\delta^*_\ell}$ for each $\hn<\theta (\ell=1,2)$ as
 $N_\ell\models
 \lq W^\hn_{\delta^*_\ell}$ is a closed unbounded subset of $\mu \rq$.
 For similar reasons $\delta^*_\ell\in W_\hn$ for each $\hn < \theta $:
 as $W_\hn\in N^{p(\hn+1)}_{\delta_\ell, \zeta}$ hence $W_\hn\in N_\ell$
 hence $W_\hn\in N_1\cap N_2$, and as $N_1, N_2\prec M^*, M^*$ has  Skolem
 functions,
 clearly $N_1\cap N_2\prec M^*$, so $W_\hn$
 is an unbounded subset of $N_1\cap N_2\cap{\mu ^+}$.
 So in $N_\ell, W_\hn$ is unbounded in
 $\delta_\ell^*=\Min [({\mu ^+}\cap N_\ell)\setminus(N_1\cap N_2)]$,
 hence $N_\ell\models\lq \delta^*_\ell\in W_\hn\rq$ hence
  $\delta^*_\ell\in W_\hn$.

  We can conclude that $\delta^*_1,\delta_2^*, \beta$
  satisfy the requirements
  (A)-(E) of $ (*)_6$ on $\delta_1, \delta_2, \xi$.
  Hence by $ (*)_7$ we have 
  $\delta_1=\delta_1^*$,
  $\delta_2=\delta^*_2$.
  But, $\zeta\in N^{\hp(\hn)}_{\delta_\ell, \zeta}\subseteq N_\ell$ hence
$\zeta< \mu \cap N_1\cap N_2$ hence $\zeta<\alpha$, so 
clause   $(*)_8$ 
contradicts the choice of
  $\zeta$, so we get a contradiction,
  thus finishing the proof of the theorem 
  \end{PROOF}   


\begin{conclusion}
\label{a20}
The condition ``have least upper bound" cannot be omitted
in\footnote{and the related works} \cite{Sh:80}.
That is:
\mn
\begin{enumerate}
\item[$\boxplus$]  there are $\bbQ$ and $\cI_\alpha(\alpha < \mu^+)$
such that:
\sn
\begin{enumerate}
\item[$(a)$]  $\bbQ$ is a forcing notion, $(< \mu)$-complete, in fact
every $\le_{\bbQ}$-increasing sequence of length $< \mu$ has a lub,
i.e. satisfies $(1)_a$
\sn
\item[$(b)$]  $\bbQ$ satisfies $(2)_b$, equivalently 
$*^1_{\mu, \mathbb{Q} }  ( b) $, see \ref{a14}
\sn
\item[$(c)$]  each $\cI_\alpha$ is a dense open subset of $\bbQ$
\sn
\item[$(d)$]  no directed $\mathbf G \subseteq \bbQ$ meets
every
$\cI_\alpha,\alpha < \mu^+$.
\end{enumerate}
\end{enumerate}
\end{conclusion}

\begin{PROOF}{\ref{a20}}
Let 
$ \kappa =2 $ and  
$\bar C$ be an $S$-club system.  If $\mu$ is a successor or just
not strongly inaccessible, choose $\bar f$ and $\bar{\cI} = \langle
\cI_\delta,\cI_{\delta,i}:\delta \in S,i < \mu\rangle$ as in
\ref{a17}, so $\bbQ = \bbQ_{(\cW,\bar f,\bar C)}$.  So $\bbQ$
satisfies clause (a) by \ref{a12}(1), satisfies clause (b) by
\ref{a14} and satisfies clauses (c),(d) by the choice of $\bar f$ and
$\bar{\cI}$.  We are left with the case $\mu$ is strongly
inaccessible, then we use \ref{a17}  the case $  yk= \mu $
  instead of  the case  $ \kappa = 2$ 
\end{PROOF}






\noindent
In \ref{a20} above we get failure when we waive
in \cite{Sh:80} the ``well met condition".
\begin{conclusion}
\label{a24}
In \ref{a20}, we may replace (a) by $(a)'$ and add (e) where:
\mn
\begin{enumerate}
\item[$(a)'$]  $\bbQ$ is a forcing notion strategically
  $(<\mu)$-complete (i.e. $(1)_c$), in fact some partial
order $\le_{\st}$ witnesses
 it in a strong way (i.e. $(1)^+_c$) ,
\sn
\item[$(e)$]  (well met)  $ (3)_a $ holds, that is 
if $p,q \in \bbQ$ are compatible then they
have a lub, (so in clause (a)' above we get $ (2)_a$.
\end{enumerate}
\end{conclusion}

\begin{PROOF}{\ref{a24}}   
We use   
a
variant of the forcing 
from  Def \ref{a7}(2)  but in clause (A)(c) there we demand  
$ h_p( \delta ) $  has a last element  (so is closed) and we
repeat the proof of 
\ref{a9}. 
Actually similarly to 
the proof of  \ref{a20},   
see \ref{c3} in particular.
In details,
this forcing notion satisfies clause
$(a)'$ by \ref{c11}(1),(2) below; clause $(b)$, i.e. $(2)_b$, by
\ref{c11}(5) below.  As for clauses (c),(d) we choose $\bar{\mathbf f}$
by \ref{a17}.  
\end{PROOF}

\begin{remark}
\label{a26}
1) In \ref{a14} and \ref{a12} we can moreover find
$\langle \cI_\varepsilon:\varepsilon
< \mu\rangle$ such that $\cI = \bigcup\limits_{\varepsilon < \mu}
\cI_\varepsilon
\subseteq 
\bbQ$ is dense and $p,q \in
\cI_\varepsilon \Rightarrow p,q$ are compatible (as in \cite{KT79}).

Why? 
Let $\cI = \{p \in \bbQ$: if $\alpha_1 < \alpha_2$ belongs to
$ v _ p $ 
then
the set $ h_p( \alpha ) $ has a last member and 
there is $\alpha \in C_{\alpha_2} \backslash \alpha_1$
such that $\otp(\alpha \cap C_{\alpha _2} 
) \in h_p(\alpha_2)\}$.  
By  \ref{a12}(2)   
we have $ \cI $  is a dense 
subset of $  \mathbb{Q} $.   

For $p \in \cI$ let
\mn
\begin{enumerate}
\item[$\bullet$]  $u_p = \{\alpha:\alpha \in v _p$ \underline{or} 
for some $\beta
\in v_p$ we have $\alpha \in C_\beta$ and $\otp(\alpha \cap C_\beta)  
\le \max(h_p(\beta))$ 
(implies 
$\otp(\alpha \cap C_\beta) 
\in h _p(\beta)$ for some $\beta \in v_p)   
\}$  
\sn
\sn
\item[$\bullet$]  $\mathbf E_1 = \{(p_1,p_2):p_1,p_2 \in \cI$ and
$\otp(u  
_{p_1})   
= \otp(u  
_{p_2})$ 
and the order preserving function $g$ 
from $u_{p_1}$ onto $u_{p_2}$ maps $v_{p_1}$ onto 
$v_{p_2},C_{ 
\alpha} \cap u_{p_1}$ onto $C_{ 
h(\alpha)} \cap 
u_{p_2}$ for $\alpha \in v_p$ and maps $h_{p_1(\alpha)}$ to 
$h_{p_2}(h(\alpha))$  or $ h_{p_2}( \alpha ) $ 
for $\alpha \in v_p\}$.  
\end{enumerate}
\mn
So $\mathbf E_1$ is an equivalence relation on $\cI$ with $\le \mu$  
classes: it is known that there is an equivalence relation $\mathbf E_2$
on $[\mu^+]^{< \mu}$ with $\mu$ equivalence classes such that $u_1 
\mathbf E_2 u_2 \Rightarrow u_1 \cap u_2 \trianglelefteq 
u_ {\ell} $.  

Easily the equivalence relation $\{(p_1,p_2):p_1 \mathbf E_1 p_2$ and
$u_{p_1} \mathbf E_2 u_{p_2}\}$ on $\cI$ is as required.  

\noindent
[Why?  Assume
$p_1 \mathbf E_2 p_2$ and $\alpha_\ell \in v_{p_\ell}$ and $\alpha_2 \in  
v_{p_2},\gamma \in C_{\alpha_1} \cap C_{\alpha_2}$ and $\otp(\gamma  
\cap C_{\alpha_\ell}) \in h_{p_\ell}(\alpha_\ell)$ for $\ell=1,2$.  
But then $\gamma \in u_{p_1} \cap u_{p_2}$ and $\gamma \in 
\dom(g_{p_1}) \cap \dom(g_{p_2})$ hence necessarily $\otp(\gamma \cap  
C_{\alpha_1}) = \otp(\gamma \cap C_{\alpha_2})$ and $g_{p_1}(\gamma) =  
g_{p_2}(\gamma)$.
Why? 
Let $ v = v_{p_1}\cup v_{p_2}$  
and choose $ \langle \gamma _\alpha : \alpha \in v\rangle $  
such that $ \gamma _ \alpha \in C_ \alpha $ and $ \delta \in v
\Rightarrow \gamma _\alpha > sup( C_\delta \cap \alpha $.  
Define $ p \in \mathbb{Q} $ by:
\begin{enumerate}
\item[$(*)_8$] 
\begin{enumerate} 
\item[(a)] $  v_p= v $ 
\item[(b)] $ u_p = u_{p_1} \cup u_{p_2} \cup \{\gamma _ \alpha : \alpha \in v \}$  
\item[(c)] $ h_p (\alpha )= h_{p_{\ell} }(\alpha ) \cup \{\gamma _\alpha  \} $ 
when $ \alpha \in v_{p_{\ell} }$
\item[(d)]  $ g_p= g_{p_1} \cup g_{p_2} \cup \{(\gamma _ \alpha , 
   \mathbf{c}_{\alpha  } \}): \alpha \in v \}   $
\end{enumerate} 
\end{enumerate} 
We can easily check that $ p $ is well defined (that is in clause (c)
if $ \alpha \in v_{p_1}\cup v_{p_2}$  then the two definitions agree; 
similarly in clause (d).
]  

\noindent
2) Note that 
for the forcing notion $ \mathbb{Q} $ from \ref{a24}, 
every 
$\le_{\bbQ}$-increasing continuous sequence of
   length $< \mu$ has a lub.
\end{remark}
\newpage

\section {Forcing axiom - non equivalence}

We use Definitions \ref{x2}, \ref{x3} freely; 
this section is dedicated to proving the following theorem:


\begin{theorem}
\label{c3}
Assume $\theta  + {\aleph_0} < \mu = \mu^{< \mu}$
and 
  $   2 \le \theta <   \mu $ 
and $\bbQ$ is adding
$\mu^+$ many 
$\mu$-Cohen.

\Then , in $\mathbf V^{\bbQ}$ we have:
\mn
\begin{enumerate}
\item[$\boxplus_{\mu,\varepsilon}$]  for some $\bbP$
\sn
\begin{enumerate}
\item[$(a)$]  $(\alpha)$   $  \quad \bbP$ is a forcing notion
\sn
\item[${{}}$] $ (\beta )$     
  $ \quad \bbP$
satisfies $(2)^\varepsilon_c$ from Definition \ref{x3} 



\sn
\item[${{}}$]  $(\gamma)$  $  \quad \bbP$ has cardinality $\mu^+$
\sn
\item[${{}}$]  $(\delta) \quad $ 
$  \bbP $ is strategically $\mu$-complete
(i.e. satisfies $(1)_{c,\mu}$ or even $(1)^+_c$),


\sn
\item[${{}}$]  $(\varepsilon) \quad$ 
we have   
$ (2)^+_{a, \mu }$



\item[${{}}$]  $(\zeta) \quad$ if $p,q \in \bbP$ are compatible \then \,
they have a lub, that is $ (3)_a$ holds; so in clause $ ( \delta ) $ we get
$ (2)_a$
\sn
\item[${{}}$]  $(\eta) \quad (2)^\varepsilon _c$  
holds 
for every limit $ \varepsilon < \mu $
\sn
\item[$(b)$]  $(\alpha)$  $  \quad \bbP$ is not equivalent to any forcing
notion satisfying $(1)_b + (2)^+_{a,=\theta}$

\hskip30pt or even just  $ (2)^ \varepsilon _{e, \theta, D}$
see Definition \ref{x2}
\sn
\item[${{}}$]  $(\beta) \quad$ moreover there is a sequence $\bar{\cI}
= \langle \cI_\alpha:\alpha < \mu^+\rangle$ of dense open

\hskip30pt  subsets of
$\bbP$ such that: if $\bbR$ is a forcing notion satisfying the

\hskip30pt conditions 
from $ (b)(\alpha ) $  
above, \then \, $\Vdash_{\bbR}$ ``there is no
directed

\hskip30pt $\mathbf G \subseteq \bbP$ which
meets $\cI_\alpha$ for $\alpha < \mu^+$".
\end{enumerate}
\end{enumerate}
\end{theorem}

\begin{remark}  \label{c4}
Hence the relevant forcing axioms are not equivalent!
\end{remark}

\begin{PROOF}{\ref{c3}}
By \ref{c11}, \ref{c22}, \ref{c25} below.

In details: let $\bar{\mathbf f}$ be from \ref{c22}(0), (i.e. after the
preliminary forcing $\bbQ$, in $\mathbf V^{\bbQ}$) and $\bbP =
\bbP_{\bar{\mathbf f}}$, 
as defined in \ref{c9}.
\smallskip

\noindent
Clause $(a)(\alpha) \quad \bbP$ a forcing notion, holds by
Definition
\ref{c9}, i.e. first statement of \ref{c11}(1).
\smallskip

\noindent
Clause $(a)(\beta)$, i.e.
for every limit ordinal $ \varepsilon < \mu $ 
the statement
$(2)^\varepsilon   
_c$ holds by \ref{c11}(7)  
\smallskip

\noindent
Clause $(a)(\gamma), ``\bbP$ of cardinality $\mu^+$",  holds by
\ref{c11}(1).
\smallskip

\noindent
Clause $(a)(\delta), (1)^+_c$ and so $\bbP$ is strategically
$\mu$-complete, by \ref{c11}(1),(2); 
\smallskip

\noindent
Clause $(a)(\varepsilon)$, means $(2)^+_ a $ 
which holds by
\ref{c11}(6)  
\smallskip

\noindent
Clause $(a)(\zeta)$, ``if $p,q$ are compatible then they have a
lub", holds by \ref{c11}(3).
\smallskip

\noindent
Clause $(a)(\eta)$ follows by $(a)(\delta), (\varepsilon)$,  

\noindent 
Clause $(b)(\alpha), ``\bbP$ not equivalent to a forcing
satisfying $(1)_b + (2)^+_{b,\theta}$" holds, by Clause $(b)(\beta)$.
\smallskip

\noindent
Clause $(b)(\beta) \quad$ the first case, 
(that is, $ {`} {`} \mathbb{R} $ satisfies $ (1)_b + (2)^+_{a,  \theta(+) })$),
  holds by \ref{c25}(2) 
because it assumption holds by \ref{c22}. 
The second  case, (that is  
$ \mathbb{R} $ satisfies  $   
(2)^\varepsilon _{e, \le \theta, D }$) 
holds by $ \ref{c25}$(3), its assumption holds by \\ref{c22}
\end{PROOF}

\begin{conclusion}
\label{c5}  
If $\theta = \cf(\theta) < \mu = \mu^{< \mu}$ \then \,
$\Ax_\mu( (1)_c + (2)^+_{a, =\theta } )$ 
does not imply
$\Ax^\theta_\mu$  and even
$\Ax  _{\mu ^{++}, \mu }((1)_c + (2)^\theta _c) $ from \ref{x5}(3). 
\end{conclusion}

\begin{PROOF}{\ref{c5}}
Let $ \lambda = \lambda ^{ \mu^+} $ 
   , $ \mathbb{Q}, \mathbb{R}  $ as in \ref{c3} and 
$ \mathbf{V} _1 = \mathbf{V} ^ \mathbb{Q} $. In $ \mathbf{V} _1 $ we can find
a
forcing notion $\bbR$ 
which 
forces $\Ax_\mu((1)c _+ (2)^+_{a, \theta(+)  })$ 
and satisfies those 
conditions, we know such $\bbR$
exists because $ (< \mu )$-support iterations
preserve the property $ (1)_c + (2)^+_{a, = \theta })$ , see \ref{z9}. 
Now $\bbR$ satisfies $(1)_b$ by
$(1)_a$ and it satisfies $(1)_{c,\theta}$ by $(1)_a$ and it satisfies
$(2)^+_{b,\theta}$ by $(2)_a + (3)_{b,\theta}$.
Now also in the universe $ \mathbf{V} _1 ^ \mathbb{R} $ the forcing
notion $ \mathbb{P} $
satisfies the conditions in $ \Ax^ \theta _ \mu $ 
from \ref{x5}.  

So by clause $ (b)( \beta ) $ of Th. \ref{c3} , in $ \mathbf{V} _1 ^ \mathbb{R} $ 
the axiom $ \Ax^\theta _ \mu $ fail as exemplified by $ \mathbb{P} $,  
so we are done proving the conclusion.
\end{PROOF}

\noindent
For this section (clearly if $\mu = \mu^{< \mu} > \aleph_0$ then there
are such objects)  
\begin{hypothesis}
\label{c6}
1) $\mu = \mu^{< \mu} > \theta  \ge 2 $ 

\noindent
2) $S = S^{\mu^+}_\mu = \{\delta < \mu^+:\cf(\delta) = \mu\}$ or $S$ just
 a stationary subset of $S^{\mu^+}_\mu$.

\noindent
3) $\bar C$ is an $S$-club sytem, see Definition \ref{a4}.

\noindent 
4) $\bar{ \mathbf{f} }$  is as in \ref{c9} but 
  $ \mathbf{f} _\delta : C_ \delta \rightarrow \theta $
\end{hypothesis}

\begin{discussion} 
1)
The main  
difference between the forcing 
in 
Def \ref{c9} below and 
the one in \ref{a7}(2)  
above is that 
\begin{enumerate}  
\item[(A)]  there the generic  gives a function $ \name{ g }$ from $ \lambda $ 
to $ \kappa $ such that for every $ \delta \in S $ for \text{`}\text{`}most" 
 $\alpha \in C_ \delta "$ we have $ \name{ g }( \alpha ) = 
 \mathbf{f}  _ \delta ( \alpha ) $
\item[(B)]  here the generic gives a function $ \name{ g }$  such that for every 
$ \delta \in S $  for 
\text{`}\text{`}most" $ \alpha \in C_ \delta $  we have
$ \mathbf{f}  _ \delta ( \alpha )\in \name{ g } ( \alpha ) $
\end{enumerate} 

\noindent 
2)xob vamud 17  ktob ktob
\end{discussion}

\begin{definition}
\label{c9}
For $\bar{\mathbf f}$ 
an $(S,\bar C,\theta  
)$-parameter, see Definition
\ref{a7}, we define a forcing notion $\bbP =
\bbP_{\bar{\mathbf f}, \theta }$ as
follows (but abusing our notation we may omit $ \theta $):
\mn
\begin{enumerate}
\item[$(A)$]  $p \in \bbP$ \Iff \, $p$ consists of (so $u_p=u$, etc.):
\sn
\begin{enumerate}
\item[$(a)$]  $u \in [\mu^+]^{< \mu}$
\sn
\item[$(b)$]  $g:u \rightarrow [\mu]^{< \theta}$, (can use $g:u
\rightarrow \theta$ because 
$\bigwedge\limits_\delta \Rang(\mathbf f_\delta)
\subseteq \theta$)
\sn
\item[$(c)$] $v \subseteq S$ of cardinality $< \mu$
\sn
\item[$(d)$]  $h$ a function with domain $v$
\sn
\item[$(e)$]  if $\delta \in v$ then
\sn
\item[${{}}$]  $(\alpha) \quad h(\delta)$ is a closed bounded
  non-empty subset of $C_\delta$
\sn
\item[${{}}$]  $(\beta) \quad h(\delta) \subseteq u$
\sn
\item[${{}}$]  $(\gamma) \quad$ if $\beta \in h(\delta)$ then
$\mathbf f_\delta(\beta) \in g(\beta)$
\end{enumerate}
\sn
\item[$(B)$]  $p \le q$, i.e. $\bbP_{\bar{\mathbf f}} \models ``p \le q"$
\Iff \,
\sn
\begin{enumerate}
\item[$(a)$]  $u_p \subseteq u_q$ and $g_p \subseteq g_q$
\sn
\item[$(b)$]  $v_p \subseteq v_q$
\sn
\item[$(c)$]  if $\delta \in v_p$ then $h_p(\delta)$ is an initial segment of
$h_q(\delta)$
\snn
\item[(d)] if $ \delta \in v_p$ and $ \alpha \in h_q ( \delta )
\setminus h_p (p)$ (hence $ h_q (\delta )  \not= h_p(\delta )$),
then   
  $ u_p \cap C_\delta \subseteq \alpha $
\end{enumerate}
\sn
\item[$(C)$]  we define $<_{\st} = <^{\bbP}_{\st}$, the strong
order by: $p <_{\st} q$ \Iff \, $p \le q$ and
\sn
\item[${{}}$]  $(d) \quad$ if $\delta \in v_p$ and $h_p(\delta) \ne
h_q(\delta)$ \then \, $\sup(h_q(\delta)) >$

\hskip25pt $\sup(\cup\{\delta \cap
C_\gamma:\gamma \in v_p \backslash \{\delta\}\})$.
\item[(D)]  Let   $ \name{ g } = \{g_p: p \in \name{ \mathbf{G} } \} $ and
$ \name{ h } = \{h_p: p \in \name{ \mathbf{G} } \} $ 
\end{enumerate}
\end{definition}

\begin{remark}  \label{c10}  
1) We may choose $ \bar{ \mathbf{f} }$  such that 
$ \mathbf{f} _\delta $ is a function to $ \mu $ or $ \kappa $
the forcing is defined similarly. It has similar properties
but it seem that the case $ \kappa = \theta $ is enough for us.

\noindent 
2)
If in clause $(A)(e)(\alpha)$ we demand only ``$h(\delta)$ is only
closed in its supremum" \then \, we get an equivalent forcing,
we lose some nice properties but gain others.
Mainly we gain in having more cases of having a lub, in
particular for increasing sequence which has an upper bound, really
any set of $< \mu$ members which has an upper bound; but we lose for
$\Delta$-systems, i.e. \ref{c11}(6).  Also we have to be more careful
in \ref{c14}.  We shall use 
the ``closed in its supremum" 
version in \S3.
\end{remark}

\begin{claim}
\label{c11}
Let $ \bar{ \mathbf{f}} $ 
be an $ (S, \bar{ C}, \theta ) $-parameter 
as in \ref{a2}, so $ $ S is a stationary 
subset of $ S^ {\mu ^+}_\mu $.

\noindent 
1) $\bbP_{\bar{\mathbf f}}$ is a forcing notion of cardinality $\mu^+$, also
$<_{\st}$ is a partial order $\subseteq <_{\bbP}$
and $p_1 \le p_2 <_{\st} p_3 \le p_4 \Rightarrow p_1 <_{\st} p_4$
and $(\forall p)(\exists q)(p <_{\st} q)$. 

\noindent
2) Any $<_{\st}$-increasing sequence in $\bbP_{\bar{\mathbf f}}$ of length $<
\mu$ has an upper bound (this is a strong/no memory version
   of strategic $\mu$-completeness), i.e. $ < _ {\st}$
exemplifies $(1)^+_c$.

\noindent
3) If $p_1,p_2 \in \bbP_{\bar{\mathbf f}}$ are compatible \then \, they
have a lub.

\noindent
4) The set 
$\{p_i:i < i(*)\}$ has a $\le$-lub in $\bbP_{\bar{\mathbf f}}$ \when
 \, $\bigwedge\limits_{i,j<i(*)}$ ($p_i,p_j$ are compatible)
and $i(*)$ is finite or $i(*) < \mu$ and for every $\delta$, the set
$\{h_{p_i}(\delta):i < i(*)$ 
satisfies  
$\delta \in v_{p_i}\}$ is finite or
at least has a maximal member.  Note this set is linearly ordered by
being an initial segment.

\noindent
4A) The set $\{p_i:i <i(*)\}$ has an ub \when \, $i(*) < \mu$ and
$\{p_i:i < i(*)\}$ is a set of pairwise compatible members of
$\bbP_{\bar{\mathbf f}}$ and $i(*)$ is finite or $i(*) < \theta$ or
at least $i(*) < \mu$ and for every limit ordinal $\alpha$ the
following set has cardinality $< \theta$:
\mn
\begin{enumerate}
\item[$\bullet$]  $\{\delta \in \bigcup\limits_{i} v_{p_i}:\alpha =
  \sup\{h_{p_i}(\delta) +1:i < i(*)$ and $\delta \in v_{p_i}\}  \} $.
\end{enumerate}
\mn
5) $\bbP_{\bar{\mathbf f}}$ satisfies clause 
$ (3)_{b, \varepsilon }$ when
$ \varepsilon < \mu $ is a limit ordinal
and $ \theta \ge 3$

\noindent
6) $\bbP_{\bar{\mathbf{f}}}$ satisfies clauses $(2)_a,(2)^+_{a, \partial} $ of 
Definition \ref{x2} 
when  
$ \partial \le \mu $. 

\noindent
7)  The forcing notion $ \mathbb{P} _{\bar{ \mathbf{f} }}$ 
satisfies $ (2)^\varepsilon_ c $ for $ \varepsilon < \mu $.
\end{claim}

\begin{PROOF}{\ref{c11}} 
1) Recall that  
$\mu = \mu^{< \mu}$ hence $\mu^+ = (\mu^+)^{< \mu}$ and
   easily $|\bbP| = \mu^+$.  Also the statements on $<_{\st}$ are
   obvious.
 What about $ \mathbb{P} _{\bar{ \mathbf{f} }}$ 
 being a quasi order?
 Assume that $ p_1 \le p_2 \le p_3$ and we shall prove
 that $ p_1 \le p_3$, clauses (a),(b),(c) of 
 \ref{c9}(B) and we shall elaborate clause (d).
 So assume  $ \delta \in v_{p_1}$ and 
 $ \alpha \in \dom (h_{p_3}(\delta )) \setminus 
 \dom (h_{p_1}(\delta )))$ and we should prove 
 that $ \dom (h_{p_1}(\delta )) \subseteq \alpha $.
 First assume $ \alpha \in \dom (h_{p_2}(\delta ))$,
 then $ p_1 \le p_2 $ implies  
  $ \dom (h_{p_1}(\delta ))  \subseteq \alpha $ 
 as required. Second assume $\alpha \notin 
 \dom (h_{p_2}(\delta ))$ then $ p_2 \le p_3$ 
 implies 
 $ \dom (h_{p_2}(\delta ))  \subseteq \alpha $
 but $ \dom (h_{p_1}(\delta )) \subseteq 
   \dom (h_{p_2}(\delta ))$ so we are done. 
 

\noindent
2) Let $\gamma < \mu$ be a limit ordinal and $\bar p = \langle p_i:i <
\gamma\rangle$ be a $<_{\st}$-increasing sequence of members of $\bbP_{\bar{\mathbf f}}$.

Let
\mn
\begin{enumerate}
\item[$(*)_1$]  $(a) \quad v_* = \bigcup\limits_{i}\{v_{p_i}:i <
\gamma\}$
\sn
\item[${{}}$]  $(b) \quad$ let $\mathbf i:v_* \rightarrow \gamma$ be
$\mathbf i(\delta) = \min\{i < \gamma:\delta \in v_{p_i}\}$
\sn
\item[${{}}$]  $(c) \quad$ let $v^*_2 = \{\delta \in
v_*$: the sequence $\langle h_{p_i}(\delta):i \in
[\mathbf i(\delta),\gamma)\rangle$ is not

\hskip25pt eventually constant$\}$
\sn
\item[${{}}$]  $(d) \quad$ for $\delta \in v^*_2$ let $\zeta_\delta =
\sup(\cup\{h_{p_i}(\delta):i \in [\mathbf i(\delta),\gamma)\}$,
\sn
\item[${{}}$]  $(e) \quad$ let $v^*_1 = v_* \backslash v^*_2$.
\end{enumerate}
\mn
We try naturally to define $p = (u_p,v_p,g_p,h_p)$ as
$\bigcup\limits_{i < \gamma} p_i$, that is
\mn
\begin{enumerate}
\item[$(*)_2$]  $(a) \quad v_p = v_* := \cup\{v_{p_i}:i < \gamma\}$
\sn
\item[${{}}$]  $(b) \quad u_p = \cup\{u_{p_i}:i < \gamma\} \cup
\{\zeta_\delta:\delta \in v^*_2\}$  
\sn
\item[${{}}$]  $(c) \quad g_p = \cup\{g_{p_i}:i < \gamma\} \cup
\{\langle \zeta_\delta,\{\mathbf f_\delta(\zeta_\delta)\}\rangle:\delta
\in v^*_2\}$
\sn
\item[${{}}$]  $(d) \quad h_p$ is a function with 
domain $v_p $ such that  
\sn
\begin{enumerate}
\item[${{}}$]  $(\alpha) \quad$ if $\delta \in v^*_1$ then
$h_p(\delta) = p_i(\delta)$ for $i < \delta$ large enough
\sn
\item[${{}}$]  $(\beta) \quad$ if $\delta
\in v^*_2$ then $h_p(\delta) = \cup\{h_{p_i}(\delta): i \in
[\mathbf i(\delta),\gamma)\} \cup \{\zeta_\delta\}$.
\end{enumerate}
\end{enumerate}
\mn
The point is to check that $p \in \bbP$, 
because 
$i < 
\gamma  
\Rightarrow p_i \le p$ 
is immediate:
\mn
\begin{enumerate}
\item[$\bullet$]  $u_p \in [\mu^+]^{< \mu}$ because $u_{p_i} \in
[\mu^+]^{< \mu}$ and $\gamma  
< \mu = \cf(\mu)$
\sn
\item[$\bullet$]  $v_p \in [S]^{< \mu}$ because $v_{p_i} \in  
[S]^{< \mu}$ and $\gamma  
< \mu = \cf(\mu)$ and $|v^*_2| \le  
\Sigma\{|v_{p_i}|:i < \gamma\} < \mu$
\sn
\item[$\bullet$]  $h_p$ is a function with domain $v_p$ such that
$\delta \in v_p   
\Rightarrow h_p(\delta)$ is a bounded closed subset of
$C_\delta$ (check the two cases)
\sn
\item[$\bullet$]  $g_p$ is a function from $u_p$ to
$ \theta $  as each
$g_{p_i}$ is a function from $u_{p_i}$  
to $\lambda$ and $\bar p$ is
$<_{\st}$-increasing and:
\sn
\item[${{}}$]  $(*) \quad$ if $\delta \in v^*_2$ then $\zeta_\delta
\notin \bigcup\limits_{i} u_{p_i}$
[Why?  This holds by \ref{c9}(B)(d) 
applied to $ p_i \le p_j $ for $ i < j < \gamma $.]
\sn
\item[${{}}$]  $(**) \quad$ if $\delta_1 \ne \delta_2 \in v^*_2$ then
$\zeta_{\delta_1} \ne \zeta_{\delta_2}$.
\end{enumerate}

[Why? see \ref{c9}(C)(e)].
\mn

\noindent
3)  
Assume $ p_1, p_2 \in \mathbb{P} $ have a common  
upper 
bound. 

\begin{enumerate} 
\item[$(*)_1$]  We define $ p \in \mathbb{P} $ as follows:
\begin{enumerate} 
\item[(a)] $ v_p = v_{p_1} \cup v_{p_2}$
\item[(b)] $ u_p = u_{p_1} \cup u_{p_2}$
\item[(c)] $ g_p = g_{p_1} \cup g_{p_2}$
\item[(d)] $ h_p $ is the function with domain $ v_p $
and for $ \delta \in v_ p $ we have
\begin{enumerate}
\item[$\bullet _1$]  if $ \delta \in v_{p_1} \setminus v_{p_2}$ 
   then $ h_p( \delta )= h_{p_1}(\delta )$
\item[$\bullet _2$] if $ \delta \in v_{p_2} \setminus v_{p_1}$ 
   then $ h_p( \delta )= h_{p_2}(\delta )$
\item[$\bullet _3$] if $ \delta \in v_{p_1} \cap v_{p_2}$ 
   then $ h_p( \delta )= h_{p_1}(\delta )\cup h_{p_2}(\delta )$
\end{enumerate} 
\end{enumerate} 
\end{enumerate}

Now
  indeed
\begin{enumerate} 
\item[$(*)_2$] $ p \in \mathbb{P} $
\end{enumerate} 

Also
\begin{enumerate} 
\item[$(*)_3$]   $ p_ {\ell}  \le p  $ for $ {\ell} = 1,2 $
\end{enumerate} 

[Why?  E.g. for clause \ref{c9}(B)(d), 
let $ \delta \in  h_p(\delta )$ and 
$ \alpha \in  h_p( \delta )
\setminus h_{p_{\ell} }(\delta )) $. 
By the choice of $ p $, necessarily $\alpha \in 
\alpha \in   h_{p_{3- {\ell} }}( \delta )
\setminus h_{p_{\ell} }(\delta ) $. Let 
$ q $ be a common upper bound of $ p_1, p_2$, 
exist by our present assumption; so clearly
$ \alpha 
\in  h_q( \delta )
\setminus h_{p_{\ell} }(\delta ) $ 
hence $ u_{p_ {\ell} } \cap C_ \delta \subseteq \alpha $ 
as promised.]

\begin{enumerate} 
\item[$(*)_4$]  if $ q $ is a common upper
bound of $  p_1, p_2 $ 
then $ p \le q $ 
\end{enumerate} 

[why? E.g. for \ref{c9}(B)(d), $ \delta \in v_p $ 
and 
$ \alpha \in 
  h_q(\delta) \setminus h_{p}(\delta ) $ 
  we should prove that 
  $ \dom (h_{p}(\delta )) \cap C_ \delta \subseteq \alpha $. 
  Now for $ {\ell} = 1,2 $ we have
  $ p_ {\ell} \le q , \delta \in v_{p_{\ell} }$ 
  and $ \alpha \in \dom (h_q (\delta )) \setminus 
  \dom ( h_{}
  p_{\ell} (\delta )) $ hence
  $ u_{p_{\ell} }\cap C_ \delta  \subseteq \alpha $. So 
  clearly 
  $$ u_p \cap C_ \delta =
  (u_{p_1}\cup u_{p_2}) =
 ( u_{p_1} \cap C_ \delta)  \cup  (u_{p_2} \cap C_ \delta)=
 \subseteq \alpha 
 $$
 
 So we are done

\noindent 
4) The proof is similar.

\noindent
4A) Similar to the proof of part (2).

\noindent
5) Easy, but we elaborate.  

So let $ p_{{\ell}, \zeta }$  for $ {\ell} \in \{ 1,2 \},
\zeta < \varepsilon $   be as in 
$ (3)_{b, \varepsilon }$. 
Let $ v = \cup \{ v_{p_{{\ell}, \zeta }}: 
  {\ell} \in \{ 1,2 \}, \zeta < \varepsilon  \} $
 
  and for $ \delta \in  v  (\subseteq S)$  let
  $ \zeta _ \delta = \sup \big( 
  \cup \{ h_{p_{{\ell}, \zeta }}(\delta ) :
  \zeta < \varepsilon , {\ell} \in \{ 1,2 \} \} $  
 and $ v_* = \{ \delta \in v:\text{ the supremum in the 
 definition of } \zeta _\delta \text{ is not obtained} \} $

Define 
\begin{enumerate} 
\item[$(*)_5$] We define $ p $ by
\begin{enumerate} 
\item[(a)]  $ u = \cup \{ u_{p_{{\ell}, \zeta }}: 
{\ell} \in \{ 1,2 \}, \zeta < \varepsilon \}  $ 
\item[(b)] $ v $ is defined above
\item[(c)]  for $ \delta \in v \setminus v_*$,  $ h(\delta )$
is defined naturally  as the eventual value
\item[(d)]  for $ \delta  \in v_*$,  $ h(\delta )$
is defined  as $ \cup \{ h_{p_{{\ell}, \zeta }}:
{\ell} \in \{ 1,2 \}, \zeta < \varepsilon  \}\cup 
\{(\zeta _ \delta , \mathbf{f} _ \delta (\zeta _ \delta )) \} $  
[so if $ \theta = 2$  there is no free choice!]
\end{enumerate} 
\end{enumerate} 

Now check.  

\noindent 
6)
For $ (2)_a $ by the proof of \ref{a14}, that is definidng 
$ \mathbf{h} $  as there, recalling part (3)

For $ (2)^\partial _{a,  \partial }  $ 
for $ \partial \le \mu $
 choose $ \mathbf{h} $ as in the proof of part (6), 
 using part (4) instead of part (3).
 
\noindent 
7)  The statement $ (2)^\varepsilon _c $  holds by 
parts (2) and (5).  

\end{PROOF}

\begin{claim}
\label{c14}
1)
$\cI_{\bar{\mathbf f},\alpha}$ is a dense open subset of
$  \mathbb{P}_ {\bar{ \mathbf{f}  } } $ 
where:
\mn
\begin{enumerate}
\item[$\bullet$]  $\cI_{\bar{\mathbf f},\alpha} = \{p \in
\bbP_{\bar{\mathbf f}}:\alpha \in u_p$ and $\alpha \in S \Rightarrow
\alpha \in v_p\}$.
\end{enumerate}

\noindent  
2) If $ \delta \in S $ and $ \alpha \in C_ \delta $ then $ \cI_{\delta, \alpha }$
is a dense open subset of $ \mathbb{P} _{\bar{ \mathbf{f} }}$ where:

\begin{enumerate}
\item[$ \bullet $] $ \cI_{\delta , \alpha } =   
\{ p  \in \mathbb{P} _{\bar{ \mathbf{f} }}:  
   \delta \in v_p $  
   and $ h_p(\delta ) \nsubseteq   \alpha \}$  
\end{enumerate} 
\end{claim}

\begin{PROOF}{\ref{c14}}  
1)
Assume $p \in \bbP_{\bar{\mathbf f}}$ and we shall find $q \in
\cI_{\bar{\mathbf f},\alpha}$ such that $p \le q$.
\medskip

\noindent
\underline{Case 1}:  If $(\alpha \notin S \vee \alpha \in v_p)$ and
$\alpha \in u_p$

Let $q=p$.
\medskip

\noindent
\underline{Case 2}: If $(\alpha \notin S \vee \alpha \in v_p)$ and
$\alpha \notin u_p$

Define $q$ by:
\mn
\begin{enumerate}
\item[$\bullet$]  $u_q = u_p \cup \{\alpha\}$
\sn
\item[$\bullet$]  $v_q = v_p$
\sn
\item[$\bullet$]  $g_q = g_p \cup \{\langle \alpha,\{0\}\rangle\}$  
\sn
\item[$\bullet$]  $h_q = h_p$.  
\end{enumerate}
\mn
Now check that   
$ q \wedge \alpha \in u_q$.  
Also $ p\le q$ is clear, e.g clause \ref{c9}(B)(d)
holds because $ \delta \in v_p \Rightarrow  h_p(\delta)
= h_q(\delta )$.  
\medskip

\noindent
\underline{Case 3}:   $\alpha \in S,\alpha \in u_p$ and $\alpha \notin v_p$

Let $\beta \in C_\alpha$ be such that $\delta \in v_p \Rightarrow
\beta > \sup(C_\delta \cap \alpha)$ and 
$ \sup (u_p \cap \delta ) < \beta $  
and 
define
$q \in \bbP_{\bar{\mathbf
f}}$ by:
\mn
\begin{enumerate}
\item[$\bullet_1$]  $u_q = u_p \cup \{\beta\}$, they may be equal, 
\sn
\item[$\bullet_2$]  $v_q = v_p \cup \{\alpha\}$
\sn
\item[$\bullet_3$]  $g_q = g_p \cup \{( \beta,\{\mathbf f_\alpha(\beta)\}
  )  \} $  
\sn
\item[$\bullet_4$]  $h_q = h_p \cup \{\langle \alpha,\{\beta\}\rangle\}$.
\end{enumerate}
\mn
Clearly $p \le q \in \cI_{\bar{\mathbf f},\alpha}$.

\noindent
2) Similarly.
\end{PROOF}

\begin{definition}
\label{c17}
We say that $\bar{\mathbf f}$ is $(\kappa,\partial)$-generic enough
\when \, $(A) \Rightarrow (B)$ and recall, 
$\bar{\mathbf f} =
\langle \mathbf f_\delta:\delta \in S  
\rangle,\mathbf
f_\delta:C_\delta \rightarrow \theta $ 
where ($\partial$ is a regular cardinality $< \mu$ 
(and
recall $\theta $ is a cardinal $ [2, \mu )$):
\mn
\begin{enumerate}
\item[$(A)$]  
\begin{enumerate} 
\item[{}] $(a) \quad E$ is a club of $\mu^+$
\sn
\item[${{}}$]  $(b) \quad \langle \alpha_{\delta,\zeta}:\zeta <
\mu\rangle$ is an increasing continuous sequence of members of
$C_\delta$

\hskip25pt  for $\delta \in E \cap S$
\sn
\item[${{}}$]  $(c) \quad h_\zeta$ is a pressing down function from $E
\cap S$ for $\zeta < \mu$
\end{enumerate} 
\sn
\item[$(B)$]  we can find $\xi < \mu$ of cofinality $\partial$ and
a sequence $\langle \delta_i:i < \kappa \rangle$ of ordinals from $E
\cap S$ such that:
\sn
\begin{enumerate} 
\item[${{}}$]  $\bullet_1 \quad$ if $\zeta < \xi$ then $h_\zeta \rest
\{\delta_i:i < \kappa\}$ is constant
\sn
\item[${{}}$]  $\bullet_2 \quad \langle
\alpha_{\delta_i,\zeta}:\zeta < \xi\rangle$ does not depend on $i <
\kappa$ hence also $\alpha = \alpha_{\delta_i,\xi}$ by continuity
\sn
\item[${{}}$]  $\bullet_3 \quad \{\mathbf f_{\delta_i}
(\alpha  
):i < \kappa\}$ 
 is equal to $ \theta $ 
where $ \alpha $ is from $ \bullet _2$.  
\end{enumerate}
\end{enumerate} 
\end{definition}

\begin{remark}
1) This is used when $*^\theta_\mu$ say: any $< \theta$ has lub inside
 the proof of \ref{c25}.

\noindent
2) For $\theta = 2$ 
as \ref{c11}(2) does not apply, we shall in \ref{c25}
need a stronger version - with the game, see
\S3.

\noindent
3) In \ref{c17} we may  
add:
\mn
\begin{enumerate}
\item[${{}}$]  $\bullet_4 \quad \{\alpha \in C_{\delta_i}:
\alpha 
  < \alpha_{\delta_i,\zeta}\}$ for some $\zeta < \xi$ does not depend on
  $i$   
\sn
\item[${{}}$]  $\bullet_5 \quad$ the $\mathbf f_{\delta_i}$'s agree on
  this  set.  
\end{enumerate}
\end{remark}

\noindent
Now in \ref{c22}, \ref{c25} 
we
shall 
arrive at  
the main point
\begin{claim}
\label{c22}
1) Assume $\bbQ$ is the forcing notion for adding $\mu^{+}$ many  
$\mu$-Cohens. \Then \, in $\mathbf V^{\bbQ}$, there is an  
$(S,\bar C,\mu)$-parameter $\bar{\bf{f}}$ which is
$(\kappa,\partial)$-generic enough  (in the sense of \ref{c17})
for our 
cardinals $ \theta \in [ 2, \mu )$ and 
regular $ \partial
\in [\aleph_0,\mu),\partial$.

\noindent
2) If $\diamondsuit_S$ then there is 
$ \bar{ \mathbf{f} } $ 
as above.
\end{claim}

\begin{PROOF}{\ref{c22}}
1) Now (modulo equivalence, so \wilog \,) $\bbQ$ can be described as follows:  
\mn
\begin{enumerate}
\item[$(*)_1$]  $(a) \quad p \in \bbQ$ \Iff \, $p$ is a function,
  $\dom(p) \in [S   
  ]^{< \mu} $ 
  and 
  for every $\delta \in
  \dom(p),p(\delta)$ is a function from some 
  strict 
  initial segment of
  $C_\delta$ 
  into $ \theta $  
  recalling $C_\delta \subseteq
\delta$ is a club of $\delta$ of order type $\mu$
\sn
\item[${{}}$]  $(b) \quad \mathbb{Q}  \models ``p \le q"$ \Iff \, $\alpha \in
  \dom(p) \Rightarrow (\alpha \in \dom(q)) \wedge (p(\alpha)
  \trianglelefteq q(\alpha))$
\sn
\item[${{}}$]  $(c) \quad$ let $\name{\mathbf f}_\delta$ for $\delta
  \in S   
  $ be $\cup\{p(\delta):p \in \name{\mathbf
    G}_{\bbQ} \text{ satisfies } \delta \in \dom(p) \}$.
\end{enumerate}
\mn
It suffices to prove $\Vdash_{\bbQ} ``\langle \name{\mathbf f}_\delta:
\delta \in S \rangle    
$ is as required".

So assume
\mn
\begin{enumerate}
\item[$(*)_2$]  $p_* \Vdash_{\mathbb{Q} }   
``\name h_\zeta$ is a pressing
  down function on $    S    
  $ for $\zeta < \mu$ and $\langle
  \name\alpha_{\delta,\zeta}:\zeta < \mu\rangle$ is increasing
  continuous sequence of members of $C_\delta$ for $\delta \in S"$.
\end{enumerate}
\mn
It suffices to find a condition $q$ above $p_*$ forcing 
that there are 
$\langle
\delta_i:i < \kappa\rangle$ and $\xi$ as in clause (B) of Definition
\ref{c17}.  For each $\delta \in S     
$ we choose
$(p_{\delta,\varepsilon},\xi_{\delta,\varepsilon},
\bar\alpha_{\delta,\varepsilon} \rangle$ by
induction on $\varepsilon < \partial$ such that:
\mn
\begin{enumerate}
\item[$(*)^3_{\delta,\varepsilon}$]  $(a) \quad p_{\delta,\varepsilon}
    \in \mathbb{Q}  
    $ is above $p_*$
\sn
\item[${{}}$]  $(b) \quad \varepsilon(1) < \varepsilon \Rightarrow
  p_{\delta,\varepsilon(1)} \le_{\bbQ} p_{\delta,\varepsilon}$
\sn
\item[${{}}$]  $(c) \quad \delta \in \dom(p_{\delta,\varepsilon})$
\sn
\item[${{}}$]  $(d) \quad \xi_{\delta,\varepsilon} =
  \otp ( \dom(p_{\delta,\varepsilon}(\delta))) $
\sn
\item[${{}}$]  $(e) \quad$ if $\varepsilon = \varepsilon(1)+1$ then
\sn
\begin{enumerate}
\item[${{}}$]  $\bullet_1 \quad p_{\delta,\varepsilon}$ forces  
a
  value $h^*_\zeta(\delta)$ to $\name h_\zeta(\delta)$ for $\zeta <
  \xi_{\delta,\varepsilon(1)}$
\sn
\item[${{}}$]  $\bullet_2 \quad p_{\delta,\varepsilon}$ forces a
value $\bar\alpha_{\delta,\varepsilon(1)}$ to $\langle
\name\alpha_{\delta,\zeta}:\zeta \le \xi_{\delta,\varepsilon(1)} +1 \rangle$
\sn
\item[${{}}$]  $\bullet_3 \quad \xi_{\delta,\varepsilon} >
\xi_{\delta,\varepsilon(1)}$ 
and 
$ \rang( \bar{ \alpha}  _{\delta ,  \varepsilon  (1)})
\subseteq \dom(p( \delta ))$.
\end{enumerate}
\end{enumerate}
\mn
There is no problem to carry the induction.  Let $\xi_\delta =
\cup\{\xi_{\delta,\varepsilon }:\varepsilon  < \partial\} <
\mu,\alpha^*_\delta =
\sup\{\dom(p_{\delta,\varepsilon}(\delta):\varepsilon <
\partial\},p_\delta = \cup\{p_{\delta,\varepsilon}:\varepsilon < \partial\}$.

Now we can define a pressing down function $h$ on $S   
$ such
that:
\mn
\begin{enumerate}
\item[$(*)_4$]   
if $\delta_1,\delta_2 \in S$     
 and
$h(\delta_1) = h(\delta_2),\varepsilon < \partial$ 
  then: 
 \begin{enumerate} 
 \item[(a)] 
  $\bar\alpha_{\delta_1,\varepsilon} =
  \bar\alpha_{\delta_2,\varepsilon}$ 
  \item[(b)] 
  for every $\alpha \in
  \Rang(\bar\alpha_{\delta_1,\varepsilon})$ we have 
  \begin{enumerate} 
  \item[$\bullet _1 $] $(C_{\delta_1}) 
  \cap \alpha = (  C_{\delta_2} \cap \alpha)$ ,
  \item[$\bullet _2$] 
  $p_{\delta_1}(\delta_1) \rest
  (C_{\delta_1} \cap \alpha) = p_{\delta_2}(\delta_2) \rest
  (C_{\delta_2} \cap \alpha)$ 
  \end{enumerate} 
  \item[(c)] 
  $h^*_\varepsilon(\delta_1) =
h^*_\varepsilon(\delta_2)$ so $\xi_{\delta_1} = \xi_{\delta_2}$ and
$p_{\delta_1, \varepsilon } \rest \delta_1 =
p_{\delta_2, \varepsilon } \rest \delta_2 $. 
   \end{enumerate} 
\end{enumerate}
\mn
Next choose an increasing sequence $\langle \delta_i:i < 
\theta \rangle$
of members of $S  
$ such that $ h $ is  
constant on $\{\delta_i:i <
\theta \}$ and $i<j \Rightarrow \dom(p_{\delta_i}) \subseteq \delta_j$.

Define $q \in \bbQ$:  
\mn
\begin{enumerate}  
\item[$(*)_5$]    
$(a) \quad \dom(q) =
  \cup\{\dom(p_{\delta_i,\varepsilon}:i < \theta ,
  \varepsilon < \theta \}$
\sn
\item[${{}}$]  $(b) \quad $ if $ i < \theta $ then 
$ q(\delta_i) =
  \cup\{p_{\delta_i,\varepsilon}(\delta_i):\varepsilon < \partial\}
  \cup \{\langle \alpha^*_\delta, i 
  \rangle\}$
  
  \item[${{}}$]  $ (c) \quad$ is $ \delta \in \dom(q) 
  \setminus \{ \delta _i :i  < \theta  \} $ 
  \then \, 
  $q(\alpha) =
  \cup\{p_{\delta_i,\varepsilon}(\alpha):\alpha \in$

\hskip25pt $\dom(p_{\delta_i,\varepsilon})\}$.
\end{enumerate}
\mn
2) Also easy.
\end{PROOF}

\begin{claim}
\label{c25}
1) There are dense sets $\cI_\alpha \subseteq \bbP =
\bbP_{\bar{\bf{f}}}$ 
for   
$ \alpha < \mu^+$, such that if $\bf{G} \subseteq
\bbP$ is directed and meets every $\cI_\alpha$, \then \, $\bf{G}$
is $\theta^+$-directed 
and even $ (< \mu )$-direccted.  

\noindent
2) If $\bar{\bf{f}}$ is $(\theta,\partial )$-generic enough and the
forcing notion $\bbR$ satisfies 
$ (1)_{c }   +(2)^{+}_{a,\theta(+)}$,
see \ref{z9}
\then \, in $\mathbf V^{\bbR}$ there is no
$ (< \mu )$-directed $\bf{G} \subseteq \bbP = \bbP_{\bar{\bf{f}}}$
meeting 
all the sets from \ref{c14}.

\noindent
3) Also there is no such $\mathbb{R}  $  
satisfying $(2)^\varepsilon_{e, \theta ,D}$
   when $\varepsilon < \mu$ is a limit ordinal 
\end{claim}

\begin{PROOF}{\ref{c25}}
1) Let $\cS = \{\bar{p}: \bar{p} \text{ is 
a directed   
sequence of conditions in } \bbP \text{ of limit length} <
\mu \} $.  
Since $\mu^{< \mu} = \mu$ and $|\bbP| = \mu^{+}$ 
it follows that   
$ |\cS| \le \mu^{+}$.
For each $\bar{p} = \langle p_i:i < 
i_*  
\rangle \in \cS$, let
$\cI_{\bar{p}} = \{q \in \bbP:q \text{ is 
either incompatible with }
p_i \text{ for some } i < i_* 
\text{ or }  p_i \le q,
\text{ for every } i <  i_* < \mu \} $.   
Since $\bbP$ is
$\mu$-strategically complete (by Claim 
    \ref{c11}(1),(2)), the set  
    $\cI_{\bar{p}}$
is dense 
  and open.  
Let $\bf{G}$ meet $\cI_{\bar{p}}$, for every $\bar{p} \in
\cS$. Then $\bf{G}$ is $\theta^{+}$-directed.
  
\noindent
2) Towards contradiction, assume $p_* \Vdash_{\bbR} ``\name{\bf{H}}
\subseteq \bbP$ is   
$ ( < \mu )$-directed,  meeting
all the sets from \ref{c14}".
Using 
 $ (1)_{c, \mu }$, fix a winning strategy
$\bf{st}$ for COM, the completeness player in the game
$ \Game_  \mu ( p^*, \mathbb{R} ) $, see Def \ref{z3}(1)
choose   
$ (E_{\zeta},\bar{q}_{\zeta},\bar{r}_{\zeta},
\bar{\mathbf{h} }_{\zeta},\bar{p}_{\zeta},  
\bar{\alpha}_{\zeta}) $ 
by induction on $ \zeta < \mu $ 
such that:   
\mn
\begin{enumerate}
\item[$(*)$]  $(a) \quad \bar{q}_{\zeta} = \langle q_{\zeta,\delta}:
\delta \in E_{\zeta} \rangle$ and $\bar{r}_{\zeta} = \langle
r_{\zeta,
\delta } 
:\delta \in E_{\zeta} \rangle$
\sn
\item[${{}}$]  $(b) \quad p_* \le q_{\zeta, \delta} \le
  r_{\zeta,\delta}$ are from $\bbR$
\sn
\item[${{}}$]  $(c) \quad \langle (q_{\xi,\delta},r_{\xi,\delta}):\xi
  \le \zeta \rangle$ is an initial segment of a play of
$\Game_{\mu }  ( p^*, \mathbb{R} ) $ 
in which

\hskip25pt   the player   
COM uses $\bf{st}$
\sn
\item[${{}}$]  $(d) \quad E_{\zeta} \subseteq \mu^{+}$ is a club
\sn
\item[${{}}$]  $(e) \quad \bf{h}_{\zeta}$ is a regressive function on
$S \cap E_{\zeta}$
\sn
\item[${{}}$]  $(f) \quad$ if $\cU \subseteq E_{\zeta} \cap S,|\cU| < \theta$
and $\bf{h}_{\zeta} \rest \cU$ is constant, then $\{r_{\zeta,\delta}:
\delta \in \cU\}$

\hskip25pt  has a lub in $\bbR$
\sn
\item[${{}}$]  $(g) \quad \bar{p}_{\zeta} = \langle p_{\zeta,\delta}:
\delta \in E_{\zeta} \rangle$
\sn
\item[${{}}$]  $(h) \quad r_{\zeta,\delta} \Vdash_{\bbR}
  ``p_{\zeta,\delta} \in \name{\mathbf H} \text{ is above }
p_{\xi,\delta} \text{ for } \xi < \zeta"$
\sn
\item[${{}}$]  $(i) \quad \bar\alpha_\zeta = \langle
  \alpha_{\delta,\zeta}:\delta \in S \cap E_{\zeta} \rangle$
\sn
\item[${{}}$]  $(j) \quad$  
$\alpha_{\delta,\zeta}$ is   a 
member of $h_{p_{\zeta,\delta}}(\delta)$
above $ \dom (h_{p_{\xi , \delta }}(\delta )$ 
for every $ \xi < \varepsilon $.
\end{enumerate}
\mn
For clauses (e)+(f), we use condition $(2)^{+}_{a, \theta}$. 

Since $\bar{\bf{f}}$ is $(\theta,\theta)$-generic enough, we 
can find 
$\langle \delta_i:i < \theta \rangle$ and $\xi$ as in Definition
\ref{c17} and let $\langle \zeta_i:i < \theta \rangle$ be
increasing with limit $\xi$.

By clause (f), for each $j < \theta$, the set $\{r_{\zeta_j,\delta_i}:
i < j\}$ has a lub $r^*_j \in \bbR$ - so necessarily $j_1 < j_2 <
\theta \Rightarrow r^*_{j_1} \le r^*_{j_2}$.  Hence the sequence
$\langle r^*_{j}:j < \theta \rangle$ has an upper bound $r_*$ (by
$(1)_{b,= \theta}$).  So $r_* \Vdash_{\bbR} \{p_{\zeta_i,\delta_j}:
i < j < \theta\} \subseteq \name{\mathbf H}$.  As $r_*
\Vdash_{\bbR} \name{\mathbf H}$ is $<\theta^{+}$-directed, we can find
some $p \in \bbP,r_{**} \ge r_*$ such that $r_{**} \Vdash_{\bbR}
p \in \name{\mathbf H}$ is an upper bound for $\{p_{\zeta_i,\delta_j}:
i < j < \theta\}$.

So, on one hand, $g_{p}(\alpha_{\delta_0, \xi})$ is a subset of $\mu$
of cardinality $< \theta$ - by the definition of $\bbP$.
On the other hand, $i < \theta \implies \alpha_{\xi,\delta_i} =
\alpha_{\xi,\delta_0}$ and $\bf{f}_{\delta_i}(\alpha_{\delta_i,\xi})
\in g_p(\alpha_{\delta_i,\xi})$.  But by Definition
\ref{c17}(B)$\bullet_3$ this is impossible.
\end{PROOF}

\begin{conclusion}
\label{c28}  
If $ \lambda = \lambda ^{< \lambda }  > 
  \mu = \mu ^{< \mu }> {\aleph_0} $ 
and 
$ \theta \not= \partial, \Delta = \cf*\partial )< \mu $ 
(and recall 
$ 2 \le \theta \le \mu $) 
\then  \, for some forcing notion $ \mathbb{R}$ we have: 
\begin{enumerate}
\item[(a)]  $ \mathbb{R} $  satisfies $ (1)_c + (2)^+_{a, = \theta }$, 
of cardinality $ \lambda $ (so add no new sequences of length $ < \mu , $
collapse no cardinality, change no cofinality and the only
possible change in
cardinal arithmetic is making $ 2^ \mu = \lambda $)
\item[(b)]  in $\mathbf{V} ^ \mathbb{R} $  we have 
$ \Ax _{\lambda, \mu } ((1)_c + (2)^+_{a, = \theta })$
\item[(c)]  in $ \mathbf{V} ^ \mathbb{R} $ the axiom 
$ \Ax ((1)_c +(2)^+_{a, \partial })$.  
fail 
\end{enumerate} 
\end{conclusion} 
\newpage

\section {Separating 
$ \Ax ^\theta_\mu,
\Ax ^\partial_\mu$ for regular
$\theta,\partial$}

Recall that $ *^ \varepsilon _{ \mu , D}$  in
\cite{Sh:546} notation is  $ (2)^ \varepsilon _{c, D}$  
here and here $ \Ax^\theta_{\mu ,D }$ is 
$\Ax (1)_c + (2)^\theta _{c,D}$, we usually  omit $ D $ 
and $ \mu $ is understood from the context.

\begin{hypothesis}
\label{d2}
1) $\mu = \mu^{< \mu}$.

\noindent
2) $S \subseteq S^{\mu^+}_\mu$ stationary.

\noindent
3) $\bar C = \langle C_\delta:\delta \in S\rangle,C_\delta$ an
   unbounded subset of $\delta$ of order type $\mu$, listed by
   $\langle \alpha^*_{\delta,\zeta}:\zeta < \mu\rangle$ in increasing
   order.

\noindent
4) $\bar{\mathbf f}$ as in \ref{d4}.

\noindent
5) $\Theta \nsubseteq \Reg \cap \mu^+$, let
  $S^{\mu^+}_\Theta = \{\delta < \mu^+:\cf(\delta) \in \Theta\}$.

\noindent
6) $2 \le \theta < \mu$ but our main interest is $\theta = 2$.
\end{hypothesis}

\begin{definition}
\label{d4}
$\bar{\mathbf f}$ is a $(\bar C,\theta )$-parameter (or uniformization problem)
when $\bar{\mathbf f} = \langle \mathbf f_\delta:\delta \in S\rangle,\mathbf
f_\delta:C_\delta \rightarrow \theta $.
\end{definition}

\begin{definition}
\label{d6}
1) We define $\bbP^1_{\bar{\mathbf f}}$ and $<_{\st}$ as in
Definition \ref{c9} but we change clause $(A)(e)$ by:
\mn
\begin{enumerate}
\item[$(e)'$]  if $\delta \in v_p$ then
\sn
\begin{enumerate}
\item[$(\alpha)$]  $h_p(\delta)$ is a bounded subset of
$C_\delta$ closed only in its supremum,
\sn
\item[$(\beta)$]  $h_p(\delta) \subseteq u_\delta$
\sn
\item[$(\gamma)$] if $\beta \in h_p(\delta)$ so
$\delta \in v_p$ \then \, $\cf(\beta) \in \Theta \Rightarrow
\mathbf f_\delta (\beta) \in g_p (\beta)$ 
(so really only

\hskip25pt $g_p \rest (u_p \cap S^{\mu^+}_\Theta)$ matters)
\sn
\item[$(\delta)$]  if $\beta \in h_p(\delta)$ and
$\cf(\beta  
) \notin S^{\mu^+}_\Theta$ then $g_p(\beta) = \emptyset$
\end{enumerate}
\sn
\end{enumerate}
\mn
2) We define $\cI^1_{\bar{\mathbf f},\alpha} \subseteq \bbP^1_{\bar{\mathbf f}}$
as in Definition \ref{c14}.
\end{definition}

\begin{claim}
\label{d9}
$\bbP^1_{\bar{\mathbf f}}$ satisfies
\mn
\begin{enumerate}
\item[$(a)$]  any increasing sequence of length $\delta < \mu,
\cf(\delta) \notin \Theta$ has a lub, i.e. $(1)_{a,=\partial}$ for
$\partial \notin \Theta$  
\sn
\item[$(b)$]   a set of pairwise compatible conditions of cardinality
$< \min(\Theta)$ has a lub - the union, i.e. $(1)_{a,
{<\min ( \Theta )}} $ 
holds.
\end{enumerate}
\end{claim}

\begin{PROOF}{\ref{d9}}
Easy.
\end{PROOF}

\begin{claim}
\label{d12}
$\bbP^1_{\bar{\mathbf f}}$ satisfies:
\mn
\begin{enumerate}
\item[$(a)$]  we have $(1)^+_c$, i.e.
\mn
\begin{enumerate}
\item[$(\alpha)$]  $<_{\st}$ is a partial order and $p_1 \le p_2
<_{\st} p_3 < p_4 \Rightarrow p_1 <_{\st} p_4$
\sn
\item[$(\beta)$]  any $<_{\st}$-increasing chain of length $< \mu$ has an ub
\end{enumerate}
\sn
\item[$(b)$]  $(\alpha) \quad$ we have $(3)_a$, i.e. if $p,q \in
\bbP^1_{\bar{\mathbf{f} }}$ are compatible then they have a lub
\sn
\item[${{}}$]  $(\beta) \quad \{p_i:i <i(*)\}$ has a lub \when \,
$i(*) < \mu$ and $\{p_i:i < i(*)\}$ is a set of

\hskip25pt  pairwise compatible
conditions and for each $\delta \in S$,

\hskip25pt  the set $\{h_{p_i}(\delta):i
<i(*)$ and $\delta \in v_{p_i}\}$ is finite; note that this set

\hskip25pt is linearly ordered by being an initial segment
\sn
\item[${{}}$]  $(\gamma) \quad \{p_i:i <i(*)\}$ has a ub \when \,
$i(*) < \mu$ and $\{p_i:i < i(*)\}$ is a set

\hskip25pt  of pairwise compatible
conditions and if $\cf(\alpha) \in \Theta$ then

\hskip25pt  $|w_{p,\alpha}| < \theta$ where
$w_{p,\alpha} = \{\delta:\delta \in \bigcup\limits_{i}
v_{p_i}$ and $\alpha = \sup\{\sup(g_{p_i}(\delta))+1$:

\hskip25pt  $i < i(*)$ and $\delta \in v_{p_i}\}\}$
\sn
\item[$(c)$]  $(\alpha) \quad (2)_a$ holds
\sn
\item[${{}}$]  $(\beta) \quad (2)^\partial_c$ that is
$*^\partial_\mu$ holds if $\partial < \mu$ is regular and $\theta \ge 2
\vee \partial \notin \Theta$
\sn
\item[$(d)$]  $(3)_{b,\varepsilon}$ holds if $\kappa =
\cf(\varepsilon) \in \mu \backslash \Theta$ 
so is regular.
\end{enumerate}
\end{claim}

\begin{PROOF}{\ref{d12}}
Like \ref{c11}, e.g.
\medskip

\noindent
\underline{Clause (a)}:  As in \ref{c11}(1),(2).
\medskip

\noindent
\underline{Clause (b)}:  Should be clear.
\medskip

\noindent
\underline{Clause (c)}:  If $\theta \ge 2$ we use $(3)_a$, i.e. the
parallel of \ref{c11}(3).  If $\theta = 1$ and $\partial \notin
\Theta$ use clause (d).
\medskip

\noindent
\underline{Clause (d)}:  Just recall $(e)(\gamma)$ of Definition \ref{d6}.
\end{PROOF}

\begin{claim}
\label{d14}
$\cI_{\bar{\mathbf \bar{ f}},\alpha}$ is a dense open subset of
$\bbP_{\bar{\mathbf f}}$ where  
\mn
\begin{enumerate}
\item[$\bullet$]  $\cI_{\bar{\mathbf f},\alpha} = \{p \in
\bbP_{\bar{\mathbf f}}:\alpha \in u_p$ and $\alpha \in S \Rightarrow
\alpha \in v_p\}$.
\end{enumerate}
\end{claim}

\begin{PROOF}{\ref{d14}}
Should be clear.
\end{PROOF}

\begin{definition}
\label{d16}
For $(\mu,\theta,\partial,D,\bar{\mathbf f})$ as in clause (A) below we define
a game $\Game_{\gn}(\bar{\mathbf f},\theta,\partial,D)$ in clause (B)
below where:
\mn
\begin{enumerate}
\item[$(A)$]  $(a) \quad \mu = \mu^{< \mu} > \partial = \cf(\partial)
  \ge \aleph_0$ and
\sn
\item[${{}}$]  $(b) \quad S \subseteq S^{\mu^+}_\mu,\bar C = \langle
  C_\delta:\delta \in S\rangle$ a club sytem
\sn
\item[${{}}$]  $(c) \quad D$ is a normal filter on $\mu^+$ to which
  $S$ belongs
\sn
\item[${{}}$]  $(d) \quad \bar{\mathbf f} = \langle \mathbf
  f_\delta:\delta \in S\rangle,\mathbf f_\delta$ is a function from
  $C_\delta$ to $\theta$
\sn
\item[$(B)$]  $(a) \quad$ a play lasts $\partial$ moves
\sn
\item[${{}}$]  $(b) \quad$ in the $\zeta$-th move, $S^\ell_\zeta \in
  D$ such that $S^2_\zeta \subseteq S^1_\zeta \subseteq S \wedge
  (\forall \xi < \zeta)(S^1_\zeta \subseteq S^2_\xi )$

\hskip25pt  and $\bar\alpha^\ell = \langle
\alpha^\ell_{\zeta,\delta}:\delta \in
  S^{\ell} _\zeta\rangle,\alpha^\ell_{\zeta,\delta} \subseteq
  C_\delta,\alpha^2_{\zeta,\delta} > \alpha^1_{\zeta,\delta} >$

\hskip25pt $\sup\{\alpha^2_{\xi,\delta}:\xi < \delta\}$
and $\mathbf h^\ell_\zeta$ pressing down functions on $S^\ell_\zeta$
\sn
\item[${{}}$]  $(c) \quad$ in the $\zeta$-th move, the anti-generic
  player chooses $S^1_\zeta,\bar\alpha^1_\zeta,\mathbf h^1_\zeta$ and
  then

\hskip25pt  the generic player  chooses  
\sn
\item[${{}}$]  $(d) \quad$ in the end of the play the generic player
  wins \when \, for some $\delta_1 < \delta_2$

\hskip25pt from $\cap\{S^2_\zeta:\zeta
  < \partial\}$ we have $\sup\{\alpha^\ell_{\zeta,\delta_1}:\zeta <
  \partial,\ell =1,2\} =$

\hskip25pt $\sup\{\alpha^\ell_{\zeta,\delta_2}:\zeta <
  \partial,\ell = 1,2\}$, call it $\alpha$ and $\mathbf
  f_{\delta_1}(\alpha) \ne \mathbf f_{\delta_2}(\alpha)$,

\hskip25pt  $\bigwedge\limits_{k < \partial} h^\ell_k(\delta_1) =
h^\ell_k(\delta_2)$.
\end{enumerate}
\end{definition}

\begin{theorem}
\label{d19}
If $\sigma \in \Theta,\theta =2$
and $\bar{\mathbf f}$ is such that in the game $\Game_{\gn}(\bar{\mathbf
  f},\theta,\sigma,D)$ from Definition \ref{d16}
the generic player wins or just does not lose,
(so $D$ a normal filter on $\mu^+,S^{\mu^+}_\mu \in D$)
\then \,:
\mn
\begin{enumerate}
\item[$(a)$]  $\bbP^1_{\bar{\mathbf f}}$ fails $\Ax^\sigma_\mu$.
\sn
\item[$(b)$]  no forcing satisfying $*^\sigma_{\mu,D}$ adds a generic to
$\bbP^1_{\bar{\mathbf f}}$, moreover
\sn
\item[$(c)$]  no forcing satisfying $*^\sigma_{\mu,D}$ adds a
$(< \mu)$-directed or just $< (\sigma^+)$-directed
$\mathbf G \subseteq \bbP^1_{\bar{\mathbf f}}$ meeting
$\mathbf{\cI}_{\bar{\mathbf f},\alpha}$ for every $\alpha < \mu^+$
(defined in \ref{c14}).
\end{enumerate}
\end{theorem}

\begin{PROOF}{\ref{d19}}
As in the proof of \ref{c25}(1), e.g.
\medskip

\noindent
\underline{Clause (c)}:

In the proof of \ref{c25}(1), we replace {\bf st} by a winning
strategy of the completeness player in the game for
$(2)^\sigma_{d,D}$, see \ref{x3} and toward contradiction assume
$\bar{\mathbf f}$ is an $(S,\bar C,\theta)$-parameter, $p_* \in
\bbP^1_{\bar{\mathbf f}}$ and $p_* \Vdash ``\name{\mathbf H} \subseteq
\bbP^1_{\bar{\mathbf f}}$ is a $(< \sigma^+)$-directed and meet every
$\cI_{\bar{\mathbf f},\alpha},\alpha < \mu^+"$.

Now for $\zeta < \sigma$ let $\mathbf Y_\zeta$ be the set of $(\bar q_\zeta,\bar
r_\zeta,\mathbf h_\zeta,E_\zeta   
,\bar p_\zeta,\bar\alpha_\zeta)$ such that:
\mn
\begin{enumerate}
\item[$\boxplus$]  $(a) \quad \langle \bar q_\xi,\bar r_\xi,\mathbf
  h_\zeta:\xi \le \zeta\rangle$ is an initial segment of a
play of the game

\hskip25pt from Definition \ref{x3} in which the
player COM uses the strategy {\bf st}
\sn
\item[${{}}$]  $(b) \quad$ so $\bar q_\zeta = \langle
  q_{\zeta,\delta}:\delta \in S_\zeta\rangle,\bar r_\zeta = \langle
  r_{\zeta,\delta}:\delta \in S_\zeta\rangle,S_\zeta \in D$ and

\hskip25pt  $S_\zeta \subseteq \{S_\xi$: for $\xi < \zeta\}$
\sn
\item[${{}}$]  $(c) \quad \bar p_\zeta = \langle
  p_{\zeta,\delta}:\delta \in S_\zeta\rangle$ and $p_{\zeta,\delta}
  \in \bbP^1_{\bar{\mathbf f}}$
\sn
\item[${{}}$]  $(d) \quad r_{\zeta,\delta} \Vdash_{\bbR}
  ``p_{\zeta,\delta} \in \name{\mathbf H}"$
\sn
\item[${{}}$]  $(e) \quad \delta \in v_{p_{\zeta ,\delta}}$
\sn
\item[${{}}$]  $(f) \quad
\sup(\dom(h_{p_{\xi,\delta}}):\xi \le \zeta\rangle$ is strictly
  increasing.
\end{enumerate}
\mn
Now we use the definition of the game $\Game_{\gn}(\bar{\mathbf
  f},\theta,\sigma,D)$ to finish as in \ref{c17}.
\end{PROOF}

\noindent
The above theorem helps for further problem as
\begin{claim}
\label{d23}
1) If a forcing notion $\bbP$ satisfies $(1)_b + (2)_a$ and $\sigma
   \in \Reg \cap \mu$ \then \, $\bbP$ satisfies $*^\sigma_\mu$,
   i.e. $(2)^\sigma_c$.

\noindent
2) If $\bbQ$ is adding $\mu^+,\mu$-Cohen $\langle
   \name\eta_\alpha:\alpha < \mu^+\rangle,\name\eta_\alpha \in {}^\mu
   \theta$ and $\theta \le \mu,\aleph_1 \le \sigma = \cf(\sigma) <
\mu,D$ a normal filter on $\mu^+$ such that $S^{\mu^+}_\mu \in D$
\then \, $\Vdash_{\bbQ} ``\langle \name\eta_\alpha:\alpha <
\mu^+\rangle$ is a $(\bar C,\mu)$-parameter and is
$(\theta,\sigma)$-generic enough and also the generic player wins in
   the game $\Game_{\gn}(\name{\bar\eta},2,\sigma,D)"$,
   pedantically replaces $D$ by the normal filter it generates.

Explain 3.9(2).
\end{claim}

\begin{conclusion}
\label{d25}
Assume $\aleph_0 \le \sigma = \cf(\sigma) < \mu = \mu^{< \mu}$ and
$\bbQ$ is the forcing notion of adding $\mu^+,\mu$-Cohens.

\noindent
1) In $\mathbf V^{\bbQ}$, there is a forcing notion $\bbP$ satisfying
   $(1)^+_c,(2)^\theta_c$ for $\theta \in \Reg \cap \mu \backslash
   \{\sigma\}$ but not $(2)^\sigma_c$.

\noindent
2) Moreover in $\mathbf V^{\bbQ}$, if $\bbR$ is a forcing notion
   satisfying $(1)_b,(2)^\sigma_c$ then it adds no generic to $\bbP$,
   in fact $|\bbP| = \mu^+$ and we should demand ``$\mathbf G \subseteq
   \bbP$ is $\sigma^+$-directed, $\mathbf G \cap \cI_\alpha \ne
   \emptyset$ for $\alpha < \mu^+$" for some dense $\cI_\alpha
   \subseteq \bbP$ for $\alpha < \mu^+$.

\noindent
3) So for some $(< \mu)$-complete $\mu^+$-c.c. forcing notion
   (satisfying $(1)_b + (2)^\sigma_c$), in $(\mathbf V^{\bbQ})^{\bbP}$
   we have $\Ax^\sigma_\mu$ but no $\mathbf G \subseteq \bbP$ as above.
\end{conclusion}

\begin{PROOF}{\ref{d25}}
In $\mathbf V^{\bbQ}$ let $\bar{\mathbf{f} } $ be from \ref{d23}(2), $\bbP$ be
$\bbP_{\bar{\mathbf f}}$ from Definition \ref{d6}.

Now (1) follows from (2).  For (2) use \ref{d19} and \ref{d9},
\ref{d12}, \ref{d14}.  For part (3) use the forcing from
\cite[1.1-1.18]{Sh:546}.
\end{PROOF}
\newpage

\bibliographystyle{amsalpha}
\bibliography{shlhetal}
\end{document}